\documentclass[11pt]{article}
\font\tenmath=msbm10
\font\sevenmath=msbm7
\font\fivemath=msbm5
\newfam\mathfam \textfont\mathfam=\tenmath
\scriptfont\mathfam=\sevenmath \scriptscriptfont\mathfam=\fivemath

\usepackage{amsmath, moreverb, hyperref}
 \hypersetup{colorlinks=true,citecolor=blue}

\usepackage{natbib, color}
\usepackage{amsfonts, fancybox}
\usepackage{amssymb}
\usepackage[american]{babel}
\usepackage{graphicx, pstricks, pst-node, color}
\usepackage{psfrag}\usepackage{subfigure}
\usepackage{flafter}
\usepackage[section]{placeins}

\newcommand{\D}{R}
\newcommand{\rank}[1]{\mathrm{rk}(#1)}

\newcommand{\DS}{\displaystyle}

\newcommand{\cA}{\mathcal{A}}

\newcommand{\cD}{\mathcal{D}}
\newcommand{\cE}{\mathcal{E}}
\newcommand{\cF}{\mathcal{F}}
\newcommand{\cG}{\mathcal{G}}
\newcommand{\cH}{\mathcal{H}}

\newcommand{\cJ}{\mathcal{J}}
\newcommand{\cK}{\mathcal{K}}
\newcommand{\cM}{\mathcal{M}}
\newcommand{\cN}{\mathcal{N}}

\newcommand{\cP}{\mathcal{P}}

\newcommand{\cX}{\mathcal{X}}

\newcommand{\cZ}{\mathcal{Z}}

\newcommand{\bX}{\mathbf{X}}
\newcommand{\bY}{\mathbf{Y}}

\newcommand{\bone}{\mathbf{1}}

\newcommand{\R}{{\rm I}\kern-0.18em{\rm R}}
\newcommand{\h}{{\rm I}\kern-0.18em{\rm H}}
\newcommand{\K}{{\rm I}\kern-0.18em{\rm K}}
\newcommand{\p}{{\rm I}\kern-0.18em{\rm P}}
\newcommand{\E}{{\rm I}\kern-0.18em{\rm E}}
\newcommand{\Z}{{\rm Z}\kern-0.18em{\rm Z}}
\newcommand{\1}{{\rm 1}\kern-0.24em{\rm I}}
\newcommand{\N}{{\rm I}\kern-0.18em{\rm N}}

\newcommand{\qq}{{\sf q}}
\newcommand{\pp}{{\sf p}}

\newcommand{\ff}{{\sf f}}

\newcommand{\pP}{{\sf P}}
\newcommand{\qQ}{{\sf Q}}

\newcommand{\pn}{\p_{\kern-0.25em n}}
\newcommand{\pnm}{\p_{\kern-0.25em n,m}}
\newcommand{\psubm}{\p_{\kern-0.25em m}}
\newcommand{\psubp}{\p_{\kern-0.25em p}}
\newcommand{\cfi}{\cF_{\kern-0.25em \infty}}

\newcommand{\symdiffsmall}{%
  \mathbin{\scriptstyle{\bigtriangleup}}}

\newcommand{\card}{\mathop{\mathrm{card}}}
\newcommand{\argmin}{\mathop{\mathrm{argmin}}}

\newcommand{\sign}{\mathop{\mathrm{sign}}}

\newcommand{\epr}{\hfill\hbox{\hskip 4pt\vrule width 5pt
                  height 6pt depth 1.5pt}\vspace{0.5cm}\par}

\newtheorem{TH1}{Theorem}[section]

\newtheorem{cor}{Corollary}[section]
\newtheorem{lem}{Lemma}[section]

\newlength{\minipagewidth}
\newcommand{\bookbox}[1]{
\par\medskip\noindent
\framebox[\textwidth]{
\begin{minipage}{\minipagewidth}
{#1}
\end{minipage} } \par\medskip }

\newcommand{\resn}{\frac{8\sigma^2}{n}\log2}

\graphicspath{{./eps/}}

\begin{document}
\title{\bf Exponential Screening and \\optimal rates of sparse
estimation}
\author{
{\sc Philippe Rigollet }
 \thanks{Princeton University.} \and {\sc Alexandre B. Tsybakov}
 \thanks{{\sc Crest} and Universit\'e Paris 6.}
}

\date{\normalsize \today}

\maketitle

\begin{abstract}
In high-dimensional linear regression, the goal pursued here is to estimate an
unknown regression function using linear combinations of a suitable
set of covariates. One of the key assumptions for the success of any
statistical procedure in this setup is to assume that the linear
combination is sparse in some sense, for example, that it involves
only few covariates. We consider a general, non necessarily linear,
regression with Gaussian noise and study a related question that is
to find a linear combination of approximating functions, which is at
the same time sparse and has small mean squared error (MSE). We
introduce a new estimation procedure, called \emph{Exponential
Screening} that shows remarkable adaptation properties. It adapts to
the linear combination that optimally balances MSE and sparsity,
whether the latter is measured in terms of the number of non-zero
entries in the combination ($\ell_0$ norm) or in terms of the global
weight of the combination ($\ell_1$ norm). The power of this
adaptation result is illustrated by showing that Exponential
Screening solves optimally and simultaneously all the problems of
aggregation in Gaussian regression that have been discussed in the
literature. Moreover, we show that the performance of the
Exponential Screening estimator cannot be improved in a minimax
sense, even if the optimal sparsity is known in advance. The
theoretical and numerical superiority of Exponential Screening
compared to state-of-the-art sparse procedures is also discussed.
\end{abstract}
\medskip

\noindent {\bf Mathematics Subject Classifications:} Primary  62G08, Secondary 62G05, 62J05, 62C20, 62G20.

\noindent {\bf Key Words:} High-dimensional regression, aggregation, adaptation, sparsity, sparsity oracle inequalities, minimax rates, Lasso, BIC.

\section{Introduction}
\setcounter{equation}{0}
\label{SEC:intro}

The theory of estimation in high-dimensional statistical models
under the sparsity scenario has been considerably developed during
the recent years. One of the main achievements was to derive {\it
sparsity oracle inequalities} (SOI), i.e., bounds on the risk of
various sparse estimation procedures in terms of the $\ell_0$ norm
(number of non-zero components) of the estimated vectors or their
approximations \citep[see][and references therein]{BicRitTsy09,
BunTsyWeg07b, BunTsyWeg07c, CanTao07, Kol08, Kol09a, Kol09b, Gee08,
ZhaHua08, Zha09}. The main message of these results was to
demonstrate that if the number of non-zero components of a
high-dimensional target vector is small, then it can be reasonably
well estimated even when the ambient dimension is larger than the
sample size. However, there was relatively few discussion of the
optimality of these bounds, mainly based on specific
counter-examples or referring to the paper by \citet{DonJoh94a},
which treats the Gaussian sequence model. The latter approach is, in
general, insufficient as we will show below. An interesting point
related to the optimality issue is that some of the bounds in the
papers mentioned above involve not only the $\ell_0$ norm but also
the $\ell_1$ norm of the target vector, which is yet another
characteristic of sparsity. Thus, a natural question is whether the
$\ell_1$ norm plays an intrinsic role in the SOI or it appears there
due to the techniques employed in the proof.

In this paper, considering the regression model with fixed design,
we will show that the role of $\ell_1$ norm is indeed intrinsic.
Once we have a ``rather general SOI" in terms the $\ell_0$ norm, a
SOI in terms of the $\ell_1$ norm follows as a consequence. This
means that we can write the resulting bound with the rate which is
equal to the minimum of the $\ell_0$ and $\ell_1$ rates (see Theorem~\ref{TH:sparseUB1}).
Unfortunately, the above mentioned ``rather general SOI" is not
available in the literature for the previously known sparse
estimation procedures. We therefore suggest a new procedure called
the {\it Exponential Screening} ({\sc es}), which satisfies the
desired bound. It is based on exponentially weighted aggregation of
least squares estimators with suitably chosen prior. The idea of
using exponentially weighted aggregation for sparse estimation is
due to \citet{DalTsy07}.  \citet{DalTsy07, DalTsy08, DalTsy09,
DalTsy10} suggested several procedures of this kind based on
continuous sparsity priors. Our approach is different because we use
a discrete prior in the spirit of earlier work by \citet{Geo86a,
Geo86b, LeuBar06,Gir08}. Unlike \citet{Geo86a, Geo86b,
LeuBar06,Gir08}, we focus on high-dimensional models and treat
explicitly the sparsity issue. Because of the high dimensionality of
the problem, we need efficient computational algorithms, and
therefore we suggest a version of the Metropolis-Hastings algorithm
to approximate our estimators (subsection~\ref{SUB:MH}). Regarding the sparsity issue, we
prove that our method benefits simultaneously from three types of
sparsity. The first one is expressed by the small rank of the design
matrix $\bX$, the second by the small number of non-zero components
of the target vector, and
the third by its small $\ell_1$ norm. 
Finally, we mention that in a work parallel to
ours, \citet{AlqLou10} consider exponentially weighted
aggregates with priors involving both discrete and continuous
components and suggest another version of the Metropolis-Hastings
algorithm to compute them.

\medskip

The contributions of this paper are the following:
\begin{itemize}
\item[(i)] We propose the {\sc es} estimator which benefits
simultaneously from the above mentioned three types of sparsity.
This follows from the oracle inequalities that we prove in
Section~\ref{SEC:UB}. We also provide an efficient and fast
algorithm to approximately compute the {\sc es} estimator and show
that it outperforms several other competitive estimators in a
simulation study.
\item[(ii)] We show that the {\sc es}
estimator attains the {\it optimal rate of sparse estimation}. To
this end, we establish a minimax lower bound which coincides with
the upper bound on the risk of the {\sc es} estimator on the
intersection of the $\ell_0$ and $\ell_1$ balls (Theorem~\ref{TH:low}).
\item[(iii)] As a consequence, we find
{\it optimal rates of aggregation} for the regression model with
fixed design. We consider the five main types of aggregation, which
are the linear, convex, model selection, subset selection and
$D$-convex aggregation, cf. \citet{Nem00, Tsy03, BunTsyWeg07c,
Lou07}. We show that the optimal rates are different from those
for the regression model with random design established in
\citet{Tsy03}. Indeed, they turn out to be moderated by the rank of
the regression matrix $\bX$. The rates are faster for the smaller
ranks. See Section~\ref{SEC:univ}.
\end{itemize}

This paper is organized as follows. After setting the problem and
the notation in Section~\ref{SEC:notation}, we introduce the {\sc
es} estimator in Section~\ref{SEC:UB} and prove that it satisfies a
SOI with a remainder term obtained as the minimum of the $\ell_0$
and the $\ell_1$ rate. This result holds with no assumption on the
design matrix ${\bf X}$, except for simple normalization. We put it
into perspective in Section~\ref{SEC:lasso_bic} where we compare it
with weaker SOI for the {\sc bic} and the Lasso estimators. In
Sections~\ref{SEC:disc} and~\ref{SEC:opt} we discuss the optimality
of SOI. In particular, Section~\ref{SEC:disc} comments why a minimax
result in \citet{DonJoh94a} with normalization depending on the
unknown parameter is not suitable to treat optimality. Instead, we
propose to consider minimax optimality on the intersection of
$\ell_0$ and $\ell_1$ balls. In Section~\ref{SEC:opt} we prove the
corresponding minimax lower bound for all estimators and show rate
optimality of the {\sc es} estimator in this sense.
Section~\ref{SEC:univ} discusses corollaries of our main results for
the problem of aggregation; we show that the {\sc es} estimator
solves simultaneously and optimally the five problems of aggregation
mentioned in (iii) above. Finally, Section~\ref{SEC:simul} presents
a simulation study demonstrating a good performance of the {\sc es}
estimator in numerical experiments.

\section{Model and notation}
\setcounter{equation}{0} \label{SEC:notation}

Let $\cZ:=\{(x_1, Y_1), \ldots, (x_n, Y_n)\}$ be a collection of
independent random couples such that $(x_i,Y_i) \in \cX \times \R$,
where $\cX$ is an arbitrary set. Assume the regression model:
$$
Y_i=\eta(x_i)+\xi_i, \ i=1, \ldots, n\,,
$$
where $\eta:\cX \to \R$ is the unknown regression function and the errors $\xi_i$ are independent Gaussian $\cN(0,\sigma^2)$.
The covariates are deterministic elements $x_1, \ldots, x_n$ of
$\cX$. Consider the equivalence relation~$\sim$ on the space of
functions $f:\cX \to \R$ such that $f\sim g$  if and only if
$f(x_i)=g(x_i)$ for all $i=1, \ldots, n$. Denote by $Q_{1:n}$ the
quotient space associated to this equivalence relation and define
the norm $\|\cdot\|$ by
$$
\|f\|^2:=\frac{1}{n}\sum_{i=1}^n f^2(x_i)\,, \quad f\in Q_{1:n}\,.
$$
Notice that $\|\cdot\|$ is a norm on the quotient space but only a
seminorm on the whole space of functions $f:\cX \to \R$. Hereafter,
we refer to it as a norm. We also define the associated inner
product
$$
\langle f, g\rangle:=\frac{1}{n}\sum_{i=1}^nf(x_i)g(x_i)\,.
$$
Let $\cH:=\{f_1, \ldots, f_M\}$, be a dictionary of $M\ge 1$ given
functions $f_j: \cX \to \R$. We approximate the regression function
$\eta$ by a linear combination ${\sf f}_\theta(x) = \sum_{j=1}^M
\theta_j f_j(x)$ with weights $\theta = (\theta_1,\dots,\theta_M)$,
where possibly $M\gg n$.

We denote by $\bX$, the $n \times M$ design matrix with elements
$\bX_{i,j}=f_j(x_i)$, $i=1, \ldots, n, j=1, \ldots, M$. We also
introduce the column vectors ${\bf
f}=(\eta(x_1),\dots,\eta(x_n))^\top$, $\bY= (Y_1,\dots,Y_n)^\top$
and $\xi= (\xi_1,\dots,\xi_n)^\top$. Let $|\cdot|_p$ denote the
$\ell_p$ norm in $\R^d$ for $p,d\ge1$ and $M(\theta)$ denote the
$\ell_0$ norm of $\theta \in \R^M$, i.e., the number of non-zero
elements of $\theta\in \R^M$. For two real numbers $a$ and $b$ we
use the notation $a\wedge b:=\min(a,b)$, $a\vee
b:=\max(a,b)$; we denote by $[a]$ the integer part of
$a$ and by $\lceil a \rceil$ the smallest integer greater than or
equal to $a$.

\section{Sparsity pattern aggregation and Exponential Screening}\label{SEC:UB}


A \emph{sparsity pattern} is a binary
vector $\pp \in \cP:=\{0,1\}^M$. The terminology comes from the fact that
the coordinates of any such vectors can be interpreted as indicators
of presence ($\pp_j=1$) or absence ($\pp_j=0$) of a given feature indexed by
$j\in\{1,\ldots,M\}$. We denote by $|\pp|$ the number of ones in the sparsity pattern $\pp$ and by $\R^\pp$ the space defined by
$$
\R^\pp=\{\theta\cdot \pp\,:\, \theta \in \R^M\}\subset \R^M\,,
$$
where $\theta \cdot \pp \in \R^M$ denotes the Hadamard product
between $\theta$ and $\pp$ and is defined as the vector $(\theta
\cdot \pp)_j=\theta_j\pp_j, j=1\ldots, M$.

For any $\pp \in \cP$, let $\hat \theta_\pp$ be any least squares
estimator defined by
\begin{equation}
\label{EQ:lse}
\begin{split}
\hat \theta_{\pp} \in \argmin_{\theta \in \R^\pp} |\bY-{\bf X}\theta|_2^2\,, 
\end{split}
\end{equation}
The following simple lemma gives an oracle inequality for the least
squares estimator.  Let
${\rm rk}(\bX)\le M\wedge n$ denote the rank of the design matrix $\bX$. 
\begin{lem}\label{LEM:lse}
Fix $\pp \in \cP$.  Then any least squares estimator $\hat \theta_\pp$
defined in \eqref{EQ:lse} satisfies
\begin{equation}
\label{EQ:pytha}
\E\|\ff_{\hat \theta_\pp}-\eta\|^2= \min_{\theta \in
\R^{\pp}}\|\ff_{ \theta}-\eta\|^2 + \sigma^2\frac{R_{\pp}}{n}
\le \min_{\theta \in \R^{\pp}}\|\ff_{ \theta}-\eta\|^2 +
\sigma^2\frac{|\pp|\wedge \D}{n}
\end{equation}
where $R_\pp$ is the dimension of the linear subspace
$\{\bX\theta\,:\, \theta \in \R^{\pp}\}$ and $\D={\rm rk}(\bX)$.
Moreover, the random variables $\xi_1, \ldots, \xi_n$ need not be
Gaussian for~\eqref{EQ:pytha} to hold.
\end{lem}
Proof of the lemma is straightforward in view of the Pythagorean
theorem.



Let $\pi=(\pi_\pp)_\pp$ be a probability measure on $\cP$, which we
will further call a prior.  The {\it sparsity pattern aggregate}
({\sc spa}) estimator is defined as $\ff_{\tilde
\theta^{\textsc{spa}}}$, where
$$
{\tilde \theta^{\textsc{spa}}}:= \frac{ \DS \sum_{\pp \in \cP} {\hat
\theta_\pp} \exp\Big(-\frac{1}{4\sigma^2} \sum_{i=1}^n
(Y_i-\ff_{\hat \theta_\pp}(x_i))^2-\frac{|\pp|}{2}\Big)\pi_{\pp}}{\DS \sum_{\pp \in
\cP}\exp\Big(-\frac{1}{4\sigma^2}\sum_{i=1}^n (Y_i -\ff_{\hat
\theta_\pp}(x_i))^2-\frac{|\pp|}{2}\Big)\pi_{\pp}}\,.
$$
As shown in \citet{LeuBar06}, the following oracle inequality holds:
\begin{equation}
\label{EQ:soiLB} \E\|\ff_{\tilde \theta^{\textsc{spa}}}-\eta\|^2 \le
\min_{\pp \in \cP: \pi_\pp\ne 0}\Big\{\E\|\ff_{\hat
\theta_\pp}-\eta\|^2 + \frac{4\sigma^2\log
(\pi_\pp^{-1})}{n}\Big\}\,.
\end{equation}
Now, we consider a specific choice of the prior $\pi$: 
\begin{equation}
\label{EQ:sprior}
\pi_\pp:=\left\{\begin{array}{ll}
\frac{1}{H}\left(\frac{|\pp|}{2eM}\right)^{|\pp|}\,,& {\rm if \ } |\pp| < \D, \\
\frac{1}{2}\,,& {\rm if \ } |\pp|=M\,,\\
0\,,& {\rm otherwise}\,,
\end{array}\right.
\end{equation}
where $\D={\rm rk}(\bX)$, we use the convention $0^0=1$ and
$H=2\sum_{k=0}^{\D} {M \choose k} \left(\frac{k}{2eM}\right)^{k} $
is a normalization factor. In this paper we study the {\sc spa}
estimator with the prior defined in~\eqref{EQ:sprior}. We call it
the {\it Exponential Screening} ({\sc es}) estimator, and denote by
$\tilde \theta^{\textsc{es}}$ the estimator $\tilde
\theta^{\textsc{spa}}$ with the prior~\eqref{EQ:sprior}. The {\sc
es} estimator is a mixture of least squares estimators corresponding
essentially to sparsity patterns $\pp$ with small size and small
residual sum of squares. Note that the weight $1/2$ is assigned to
the least squares estimator on the whole space (case where
$|\pp|=M$) and can be changed to any other constant in $(0,1)$
without modifying the rates presented below, as long as $H$ is
modified accordingly.

Since $ {M \choose k} \le \big(\frac{eM}{k}\big)^{k}$, we obtain
that $H \le 4$.
Using this and considering separately the cases
$|\pp|\le 1$ and $|\pp|\ge 2$, we obtain that the remainder term
in~\eqref{EQ:soiLB} satisfies
\begin{eqnarray}\label{EQ:remainder}
\frac{4\sigma^2\log (\pi_\pp^{-1})}{n}&\le& \frac{4\sigma^2}{n}
\left[ |\pp|\log\Big(\frac{2eM}{ |\pp|\vee 1}\Big) +\log 4\right]\\
&\le& \frac{8\sigma^2|\pp|}{n}\log\Big(1+\frac{eM}{|\pp|\vee 1}\Big)
+\frac{8\sigma^2}{n}\log2\,\nonumber
\end{eqnarray}
for sparsity patterns $\pp$ such that $ |\pp|< \D$. Together
with~\eqref{EQ:soiLB}, this inequality yields the following theorem.

\begin{TH1}\label{TH:sparseUB}
For any $M \ge 1, n\ge 1$, the Exponential Screening estimator
satisfies the following sparsity oracle inequality
\begin{equation}
\label{EQ:soiTH} \E\|\ff_{\tilde \theta^{\textsc{es}}}-\eta\|^2 \le
\min_{\theta \in \R^{M}}\left\{\|\ff_{ \theta}-\eta\|^2
+\frac{\sigma^2 R}{n}\wedge  \frac{9\sigma^2
M(\theta)}{n}\log\left(1+\frac{eM}{ M(\theta)\vee 1}\right)\right\}+
\frac{8\sigma^2}{n}\log2
\end{equation}
 where $\D\le
M\wedge n$ denotes the rank of the design matrix $\bX$.
\end{TH1}

\noindent {\sc Proof.} {
Combining the result of Lemma~\ref{LEM:lse} and ~\eqref{EQ:soiLB}
with the sparsity prior defined in~\eqref{EQ:sprior}, we obtain that
$\E\|\ff_{\tilde \theta^{\textsc{es}}}-\eta\|^2$ is bounded from
above by
\begin{equation}\label{EQ:minRd}
 \min_{\substack{\theta \in \R^{M}\\M(\theta)< \D}}
 \left\{\|\ff_{\theta}-\eta\|^2 +
9\sigma^2\frac{
M(\theta)}{n}\log\left(1+\frac{eM}{ M(\theta)\vee
1}\right)\right\}+ \frac{8\sigma^2}{n}\log2
 \,,
\end{equation}
and by
\begin{equation}
\label{EQ:minRd2}
 \min_{\substack{\theta \in \R^{M}}}
 \left\{\|\ff_{\theta}-\eta\|^2 +\sigma^2\frac{R}{n}\right\}+ \frac{4\sigma^2}{n}\log2
 \,.
\end{equation}
Combining~\eqref{EQ:minRd} and~\eqref{EQ:minRd2} concludes the proof.
%
\epr

}

%

An interesting corollary of Theorem~\ref{TH:sparseUB} is obtained
for the linear regression model where it is assumed that
$\eta=\ff_{\theta^*}$ for some $\theta^* \in \R^M$. In this
case~\eqref{EQ:soiTH} yields
$$
\E\|\ff_{\tilde \theta^{\textsc{es}}} -\ff_{\theta^*}\|^2
\le\frac{\sigma^2 R}{n}\wedge  \frac{9\sigma^2
M(\theta^*)}{n}\log\left(1+\frac{eM}{ M(\theta^*)\vee 1}\right)+
\frac{8\sigma^2}{n}\log2\,.
$$
However, even in this parametric case, Theorem~\ref{TH:sparseUB} provides a stronger result. Indeed, if there exists $\theta' \in \R^M$, such that
\begin{equation}
\label{EQ:tradeoff} \|\ff_{ \theta'}-\ff_{\theta^*}\|^2 +
\frac{\sigma^2 R}{n}\wedge  \frac{9\sigma^2
M(\theta')}{n}\log\left(1+\frac{eM}{ M(\theta')\vee 1}\right)
 < \frac{\sigma^2 R}{n}\wedge  \frac{9\sigma^2 M(\theta^*)}{n}\log\left(1+\frac{eM}{ M(\theta^*)\vee
1}\right)\,,
\end{equation}
then Theorem~\ref{TH:sparseUB} gives a tighter bound on
$\E\|\ff_{\tilde \theta^{\textsc{es}}} -\ff_{\theta^*}\|^2$. A
vector $\theta' \in \R^M$ that satisfies~\eqref{EQ:tradeoff} exists
when $\ff_{\theta^*}$ can be well approximated by $\ff_{\theta'}$
and $\theta'$ is much sparser than $\theta^*$.

While the sparsity oracle inequality~\eqref{EQ:soiTH} indicates that
the {\sc es} estimator adapts to the underlying sparsity when
measured in terms of the number of non-zero coefficients
$M(\theta)$, it is also adaptive to the sparsity when measured in
terms of the $\ell_1$ norm $|\theta|_1=\sum_j |\theta|_j$. This can
become an advantage when $\theta$ has many small coefficients so
that $|\theta|_1\ll M(\theta)$. Indeed, the following theorem shows
that the {\sc es} estimator also enjoys adaptation in terms of its
$\ell_1$ norm.
\begin{TH1}
\label{TH:sparseUB1}  Assume that $\max_{1\le j\le M}\|f_j\|\le 1$.
Then for any $M \ge 1, n\ge 1$ the Exponential Screening estimator
satisfies
\begin{equation}
\label{EQ:soiCOR} \E\|\ff_{\tilde \theta^{\textsc{es}}}-\eta\|^2
\le\min_{\theta \in \R^{M}}\left\{\|\ff_{ \theta}-\eta\|^2 + \bar
\varphi_{n,M}(\theta)\right\}+\frac{\sigma^2}{n}(9\log(1+eM)+8\log2).
\end{equation}
where $\varphi_{n,M}(0):=0$ and, for $\theta\neq 0$,
\begin{equation}
\label{EQ:phi} \varphi_{n,M}(\theta):= \frac{\sigma^2 R}{n} \wedge \frac{9\sigma^2
M(\theta)}{n} \log\left(1+\frac{eM}{ M(\theta) \vee
1}\right)\wedge\frac{11\sigma|\theta|_1}{\sqrt{n}} \sqrt{\log
\left(1+\frac{3eM\sigma}{|\theta|_1 \sqrt{n}}\right)} \,.
\end{equation}
Furthermore, for any $\theta\in\R^M$, such that $\langle \ff_\theta,
\eta \rangle \le \|\ff_\theta\|^2$, we have
\begin{equation}
\label{EQ:phi1} \E\|\ff_{\tilde \theta^{\textsc{es}}}-\eta\|^2
\le\|\ff_{ \theta}-\eta\|^2 + \psi_{n,M}(\theta) +\resn\,
\end{equation}
where $\psi_{n,M}(0):=0$ and, for $\theta\neq 0$,
\begin{equation}
\label{EQ:psi}
 \psi_{n,M}(\theta):=\frac{\sigma^2 R}{n} \wedge \frac{
9\sigma^2 M(\theta)}{n} \log\left(1+\frac{eM}{
M(\theta) \vee 1}\right) \wedge  \frac{11\sigma |\theta|_1}{\sqrt{n}}
\sqrt{\log \left(1+\frac{3eM\sigma}{|\theta|_1 \sqrt{n}}\right)}\wedge
4|\theta|_1^2\,.
\end{equation}
In particular, if there exists $\theta^* \in \R^M$ such that
$\eta=\ff_{\theta^*}$, we have
\begin{equation}
\label{EQ:soiCORparma} \E\|\ff_{\tilde
\theta^{\textsc{es}}}-\ff_{\theta^*}\|^2 \le \psi_{n,M}(\theta^*)
+\resn \,.
\end{equation}
\end{TH1}
The proof of Theorem~\ref{TH:sparseUB1} is obtained by combining
Theorem~\ref{TH:sparseUB} and Lemma~\ref{LEM:maurey} in the
appendix. For brevity, the constants derived from
Lemma~\ref{LEM:maurey} are rounded up to the closest integer.


It is easy to see that in fact Lemma~\ref{LEM:maurey} implies a more
general result. Not necessarily $\ff_{\tilde \theta^{\textsc{es}}}$
but, in general, any estimator satisfying a SOI of the
type~\eqref{EQ:soiTH} also obeys the oracle inequality of the
form~\eqref{EQ:soiCOR}, i.e., enjoys adaptation simultaneously in
terms of the $\ell_0$ and $\ell_1$ norms. This remains still a
theoretical proposal, since we are not aware of estimators
satisfying~\eqref{EQ:soiTH} apart from $\ff_{\tilde
\theta^{\textsc{es}}}$. However, there are estimators for which
coarser versions of~\eqref{EQ:soiTH} are available as discussed in
the next section.


\section{Sparsity oracle inequalities for the BIC and Lasso estimators}
\label{SEC:lasso_bic} The aim of this section is to put
Theorem~\ref{TH:sparseUB1} in perspective by discussing weaker
results in the same spirit for two popular estimators, namely, the
BIC and the Lasso estimators.

We consider the following version of the BIC estimator, cf.
\citet{BunTsyWeg07c}:
\begin{equation}
\label{EQ:BIC}
\begin{split}
\hat \theta^{\textsc{bic}} \in \argmin_{\theta \in \R^M}
\left\{\frac1{n}|\bY-{\bf X}\theta|_2^2+{\rm pen}(\theta)\right\}\,,
\end{split}
\end{equation}
where
$$
\text{pen}(\theta) := \frac{2 \sigma^2}{n} \left\{
1+\frac{2+a}{1+a}\sqrt{ L(\theta)} + \frac{1+a}{a} L(\theta)
\right\}  M(\theta)\,,
$$
with for some $a>0$ and $L(\theta)=2\log\left(\frac{eM}{
M(\theta)\vee 1}\right)$. Combining Theorem~3.1
in~\citet{BunTsyWeg07c} and Lemma~\ref{LEM:maurey} in the appendix
we get the following corollary.
\begin{cor}
\label{COR:BIC} Assume that $\max_{1\le j\le M}\|f_j\|\le 1$. Then
there exists a positive numerical constant $C$ such that for any $M
\ge 2, n\ge 1$ and any $a>0$ the {\sc bic} estimator satisfies
\begin{equation}
 \E\|\ff_{\hat \theta^{\textsc{bic}}}-\eta\|^2
\le (1+a)\min_{\theta \in \R^{M}}\left\{\|\ff_{ \theta}-\eta\|^2 +
C\frac{1+a}{a} \varphi_{n,M}(\theta)\right\}\,+
\frac{C\sigma^2}{n}\,, \label{EQ:corBIC}
\end{equation}
where $\varphi_{n,M}$ is defined in~\eqref{EQ:phi}.
\end{cor}
We note that Theorem~3.1 in~\citet{BunTsyWeg07c} is stated with
$R=M$ and with the
additional assumption that all the functions $f_j$ are uniformly
bounded. Nevertheless, this last condition is not used in the proof
in~\citet{BunTsyWeg07c}, and the result trivially extends to the
framework that we consider here. The SOI~\eqref{EQ:corBIC} ensures
adaptation to sparsity simultaneously in terms of the $\ell_0$ and
$\ell_1$ norms. However, it is less precise than the SOI in
Theorem~\ref{TH:sparseUB1} because the leading constant $(1+a)$ is
strictly greater than 1 and the rate deteriorates as the leading
constant approaches 1, i.e., as $a\to 0$. Also the computation of
the {\sc bic} estimator is a hard combinatorial problem, exponential
in $M$, and it can be efficiently solved only when the dimension $M$
is small.

Consider now the Lasso estimator $\hat \theta^{\textsc{l}}$, i.e., a solution
of the minimization problem
\begin{equation}
\label{EQ:lasso}
\begin{split}
\hat \theta^{\textsc{l}} \in \argmin_{\theta \in \R^M} \left\{
\frac1{n}|\bY-{\bf X}\theta|_2^2 + \lambda |\theta|_1\right\}\,,
\end{split}
\end{equation}
where $\lambda>0$ is a tuning parameter. This problem is convex, and
there exist several efficient algorithms of computing $\hat
\theta^{\textsc{l}}$ in polynomial time.

Our aim here is to present results in the spirit of
Theorem~\ref{TH:sparseUB1} for the Lasso. They have a weaker form
than for the {\sc es} estimator and for the {\sc bic}. In the next
theorem, we give a SOI in terms of the $\ell_1$ norm that is similar
to those that we have presented for the {\sc es} and {\sc bic}
estimators but it is stated in probability rather than in
expectation and the logarithmic factor in the rate is less accurate.
Note that it does not require any restrictive condition on the
dictionary $f_1,\dots,f_M$.
\begin{TH1}
\label{TH:lasso1} Assume that $\max_{1\le j\le M}\|f_j\|\le 1$. Let
$M\ge 2, n\ge1$ and let $\hat \theta^{\textsc{l}}$ be the Lasso estimator
defined by~\eqref{EQ:lasso} with $\lambda=A\sigma\sqrt{\frac{\log
M}{n}}$, where $A>2\sqrt{2}$. Then with probability at least
$1-M^{1-A^2/8}$ we have
\begin{equation}
\label{EQ:soi:lasso1} \|\ff_{\hat \theta^{\textsc{l}}}-\eta\|^2 \le
\min_{\theta \in \R^{M}}\left\{\|\ff_{ \theta}-\eta\|^2 +
2A\sigma\frac{|\theta|_1}{\sqrt{n}} \sqrt{\log M}\right\}\,.
\end{equation}
\end{TH1}
{\sc Proof.} From the definition of $\hat \theta^{\textsc{l}}$ by a
simple algebra we get
$$
\|\ff_{\hat \theta^{\textsc{l}}}-\eta\|^2 \le \|\ff_{
\theta}-\eta\|^2 + \frac{2}{n}\left|(\hat
\theta^{\textsc{l}}-\theta)^\top\bX^\top\xi\right| +
\lambda\big(|\theta|_1 - |\hat \theta^{\textsc{l}}|_1\big), \quad
\forall \ \theta \in \R^{M}.
$$
Next, note that $P(\mathcal{A}) \ge 1-M^{1-A^2/8}$ for the random
event $\mathcal{A}= \Big\{ \big|\frac{2}{n}\bX^\top\xi\big|_\infty
\le \lambda\Big\}$ \citep[cf.][eq. (B.4)]{BicRitTsy09}. Therefore,
$$
\|\ff_{\hat \theta^{\textsc{l}}}-\eta\|^2 \le \|\ff_{ \theta}-\eta\|^2 +
\lambda|\hat \theta^{\textsc{l}}-\theta|_1 + \lambda\big(|\theta|_1 - |\hat
\theta^{\textsc{l}}|_1\big), \quad \forall \ \theta \in \R^{M},
$$
with probability at least $1-M^{1-A^2/8}$. Thus,
\eqref{EQ:soi:lasso1} follows by the triangle inequality and the
definition of $\lambda$. \epr

The rate $\frac{|\theta|_1}{\sqrt{n}} \sqrt{\log M}$
in~\eqref{EQ:soi:lasso1} is slightly worse than the corresponding
$\ell_1$ term of the rate of {\sc es} estimator, cf.~\eqref{EQ:phi}
and~\eqref{EQ:psi}.

In contrast to Theorem~\ref{TH:lasso1}, a SOI in terms of the
$\ell_0$ norm for the Lasso is available only under strong
conditions on the dictionary $f_1,\dots,f_M$. Following
\citet{BicRitTsy09}, we say that the restricted eigenvalue condition
RE($s$,$c_0$) is satisfied for some integer $s$ such that $1\le s\le
M$,  and a positive number $c_0$ if we have:
$$
\kappa(s, c_0) :=\min_{\substack{J_0\subseteq \{1,\dots,M\},\\
\\ |J_0|\le s}} \ \ \min_{\substack{\Delta\neq 0, \\ |\Delta_{J_0^c}|_1\le c_0|\Delta_{J_0}|_1}}
\ \ \frac{|\bX\Delta|_2}{\sqrt{n}|\Delta_{J_0}|_2} \,> \, 0.
$$
Here $|J|$ is the cardinality of the index set $J$ and we denote by
$\Delta_J$ the vector in $\R^M$ that has the same coordinates as
$\Delta$ on $J$ and zero coordinates on the complement $J^c$ of $J$.
A typical SOI in terms of the $\ell_0$ norm for the Lasso is given
in Theorem~6.1 of~\citet{BicRitTsy09}. It guarantees that, under the
condition RE($s$,$3+4/a$) and the assumptions of
Theorem~\ref{TH:lasso1}, with probability at least $1-M^{1-A^2/8}$,
we have
\begin{equation}
\label{EQ:soi:lasso2} \|\ff_{\hat \theta^{\textsc{l}}}-\eta\|^2 \le
(1+a)\min_{\theta \in \R^{M}: M(\theta)\le s}\left\{\|\ff_{
\theta}-\eta\|^2 + \frac{C(1+a)}{a
\kappa^2(s,3+4/a)}\frac{M(\theta)\log M}{n} \right\}\,,
\end{equation}
for all $a>0$ and some constant $C>0$ depending only on $A$ and
$\sigma$. This oracle inequality is substantially weaker than
\eqref{EQ:soiCOR} and \eqref{EQ:corBIC}. Indeed, it is valid under
assumption RE($s$,$3+4/a$), which is a strong condition.
Furthermore, the rank of the matrix $\bX$ does not appear, the minimum in \eqref{EQ:soi:lasso2} is taken
over the set of sparsity $s$ linked to the properties of the matrix
$\bX$, and the minimal restricted eigenvalue $\kappa(s,3+4/a)$
appears in the denominator. This contrasts with inequalities
\eqref{EQ:soiCOR}, \eqref{EQ:corBIC} and \eqref{EQ:soi:lasso1} which
hold under no assumption on $\bX$, except for simple normalization:
$\max_{1\le j\le M}\|f_j\|\le 1$. Finally, the leading constant in
\eqref{EQ:soi:lasso2} is strictly larger than $1$, and the same
comments as for the {\sc
bic} apply in this respect.  

\section{Discussion of the optimality}

\setcounter{equation}{0}

\subsection{Deficiency of the approach based on function normalization}
\label{SEC:disc}

Section \ref{SEC:UB} provides upper bounds on the risk of {\sc es}
estimator. A natural question is whether these bounds are optimal.
At first sight, to show the optimality it seems sufficient to prove
that there exists $\theta\in \R^M$ and $\eta$ such that, for any
estimator $T$,
\begin{equation*}
\E\|T-\eta\|^2 \ge \|\ff_{ \theta}-\eta\|^2 + c\psi_{n,M}(\theta)\,,
\end{equation*}
where $c>0$ is some constant independent of $n$ and $M$. This can be
also written in the form
\begin{eqnarray}\label{low:1}
&&\inf_{T}\sup_{ \eta}\sup_{\theta\in \R^M}
\frac{\E\|T-\eta\|^2-\|\ff_{ \theta}-\eta\|^2
}{\psi_{n,M}(\theta)}\ge c
\end{eqnarray}
where $\inf_{T}$ denotes the infimum over all estimators. We note
that it is possible to prove~\eqref{low:1} under some assumptions on
the dictionary $f_1,...,f_M$. However, we do not consider this type
of results because they do not lead a valid notion of optimality.
Indeed, since the rate $\psi_{n,M}(\theta)$ is a function of
parameter $\theta$, there exists infinitely many different rate
functions $\psi_{n,M}(\cdot)$ for which~\eqref{low:1} can be proved
and complemented by the corresponding upper bounds. To illustrate
this point, consider a basic example defined by the following
conditions:\\
(i) $M=n$,\\ (ii) $\eta=\ff_{\theta^*}$ for some $\theta^*\in
\R^n$,\\ (iii) the Gram matrix $\Psi=\bX^\top\bX/n$ is equal to the
$n\times n$ identity matrix,\\ (iv) $\sigma^2=1$.\\
This will be
further referred to as the {\it diagonal model}. It can be
equivalently written as a Gaussian sequence model
\begin{equation}\label{EQ:gauss_seq_mod}
y_i=\theta_i + \frac{1}{\sqrt n}\varepsilon_i, \quad i=1,\dots,n,
\end{equation}
where $(y_1,\dots,y_n)^\top = \bX^\top \bY/n$ and
$\varepsilon_1,\dots,\varepsilon_n$ are i.i.d. standard Gaussian
random variables.

Clearly, estimation of $\eta$ in the diagonal model is equivalent to
estimation of $\theta^*$ in model~\eqref{EQ:gauss_seq_mod}, and we
have the isometry $\|\ff_\theta - \eta\|= |\theta-\theta^*|_2$. Moreover, it is easy to see
that we can consider w.l.o.g. only estimators $T$ of the form
$T=\ff_{\hat\theta}$ for some statistic $\hat\theta$, and that~\eqref{low:1} for the diagonal model follows from a simplified bound
\begin{eqnarray}\label{low:2}
&&\inf_{\hat\theta}\sup_{\theta\in \R^n}
\frac{E_{\theta}|\hat\theta-\theta|_2^2 }{\psi_{n}(\theta)}\ge c,
\end{eqnarray}
where we write $E_{\theta}$ to specify the dependence of the
expectation upon $\theta$, $\inf_{\hat\theta}$ denotes the infimum
over all estimators, and for brevity
$\psi_{n}(\theta)=\psi_{n,n}(\theta)$.

Results of the type~\eqref{low:2} are available in~\citet{DonJoh94a}
where it is proved that, for the diagonal model,
\begin{eqnarray}\label{low:3}
&&\inf_{\hat\theta}\sup_{\theta\in \R^n}\
\frac{E_{\theta}|\hat\theta-\theta|_2^2 }{ \psi_{n}^{01}(\theta)}
=1+o(1),
\end{eqnarray}
as $n\to\infty$, where
\begin{equation}\label{low:4}
\psi_{n}^{01}(\theta) = 2\log n\bigg\{ \frac{1}{n} +
\sum_{i=1}^n\min\bigg(\theta_i^2, \frac{1}{n}\bigg) \bigg\}\,.
\end{equation}
The expression in curly brackets in~\eqref{low:4} is
 the risk of 0-1 (or ``keep-or-kill") oracle, i.e., the minimal
 risk of the estimators $\hat\theta$
whose components $\hat\theta_j$ are either equal to $y_j$ or to 0. A
relation similar to~\eqref{low:3}, with the infimum taken over a
class of thresholding rules, is proved in~\citet{FosGeo94}.

The result~\eqref{low:3} is often wrongly interpreted as the fact
that the factor $2\log n$ is the ``unavoidable" price to pay for
sparse estimation. In reality this is not true, and~\eqref{low:3}
cannot be considered as a basis of valid notion of optimality.
Indeed, using the results of Section \ref{SEC:UB}, we are going to
construct an estimator whose risk is $O(\psi_{n}^{01}(\theta))$ for
all $\theta$, and is of order $o(\psi_{n}^{01}(\theta))$ for some
$\theta$, cf. Theorem~\ref{TH:DonJohn1} below. So, this estimator
improves upon~\eqref{low:3} not only in constants but in the rate;
in particular, the exact asymptotic constant appearing
in~\eqref{low:3} is of no importance. The reason is that the lower
bound for~\eqref{low:3} in~\citet{DonJoh94a} is proved by
restricting $\theta$ to a small subset of $\R^n$, and the behavior
of the risk on other subsets of $\R^n$ can be much better. 

Define the rate
\begin{equation*}
\label{EQ:rate} \psi_{n}^*(\theta)=\min\left[\frac{M(\theta)\log
n}{n}, |\theta|_1\sqrt{\frac{\log n}{n}}, |\theta|_1^2\right]+
\frac1{n}\,,
\end{equation*}
which is an asymptotic upper bound on the rate in~\eqref{EQ:psi} for
$M=n$, $n\to\infty$.


\begin{TH1}
\label{TH:DonJohn} Consider the diagonal model.
 Then the Exponential Screening estimator satisfies
\begin{eqnarray}\label{EQ:THDonJohn1}
&&\limsup_{n\to\infty}\sup_{\theta\in \R^n} \frac{E_{\theta}|\tilde
\theta^{\textsc{es}}-\theta|_2^2 }{\psi_{n}^{01}(\theta)}\le 2,
\end{eqnarray}
and
\begin{eqnarray}\label{EQ:THDonJohn3}
&&\liminf_{n\to\infty}\inf_{\theta\in\R^n} \frac{E_{\theta}|\tilde
\theta^{\textsc{es}}-\theta|_2^2}{\psi_{n}^{01}(\theta)} = 0.
\end{eqnarray}
Furthermore,
\begin{eqnarray}\label{EQ:THDonJohn2}
&&\lim_{n\to\infty}\inf_{\theta\in \R^n}
\frac{\psi_{n}^*(\theta)}{\psi_{n}^{01}(\theta)} = 0.
\end{eqnarray}
\end{TH1}
{\sc Proof.} We first prove~\eqref{EQ:THDonJohn1}.
From~\eqref{EQ:soiLB}, Lemma~\ref{LEM:lse} and~\eqref{EQ:remainder}
we obtain
\begin{equation*}
E_{\theta}|\tilde \theta^{\textsc{es}}-\theta^*|_2^2 \le
\min_{\theta \in \R^n}\Big\{ |\theta-\theta^*|_2^2 +
\frac{M(\theta)}{n}(1 + 4\log (2en))\Big\} + \frac{4\log2}{n} \,.
\end{equation*}
for any $\theta^*\in\R^n$. Let $\bar \theta\in\R^n$ be the vector
with components $\bar \theta_j=\theta^*_jI(|\theta^*_j|>1/\sqrt{n})$
where $I(\cdot)$ denotes the indicator function. Then
$$
|\bar\theta-\theta^*|_2^2 = \sum_{j=1}^n |\theta^*_j|^2
I(|\theta^*_j|\le 1/\sqrt{n}),
$$
and
$$
\frac{M(\bar\theta)}{n} = \sum_{j=1}^n \frac1{n} I(|\theta^*_j|>
1/\sqrt{n}).
$$
Therefore,
\begin{equation*}
\label{EQ:proofDonJon} E_{\theta}|\tilde
\theta^{\textsc{es}}-\theta^*|_2^2 \le (1 + 4\log (2en))
\sum_{j=1}^n \min\left(|\theta^*_j|^2, \frac1{n}\right) +
\frac{4\log 2}{n} \,,
\end{equation*}
which implies~\eqref{EQ:THDonJohn1}. Next,~\eqref{EQ:THDonJohn3} is
an immediate consequence of~\eqref{EQ:THDonJohn2}. To
prove~\eqref{EQ:THDonJohn2} we consider, for example, the set
$\Theta_n = \Big\{\theta\in\R^n: \, a/\sqrt{n} \le |\theta_j| \le
b/\sqrt{n} \ \text{\rm for all} \ \theta_j\ne 0\Big\}$ where
$0<a<b<\infty$ are constants. For all $\theta\in\Theta_n$ we have
$$
\psi_{n}^*(\theta) \le |\theta|_1\sqrt{\frac{\log n}{n}} +
\frac1{n}\le b M(\theta) \frac{\sqrt{\log n}}{n}+\frac1{n}\,,
$$
and
$$
\psi_{n}^{01}(\theta) \ge 2(\min(a^2,1)M(\theta)+1) \frac{\log
n}{n}\,,
$$
so that
\begin{equation}\label{EQ:infty}
\lim_{n\to\infty}\sup_{\theta\in \Theta_n}
\frac{\psi_{n}^*(\theta)}{\psi_{n}^{01}(\theta)} = 0.
\end{equation}
Hence, \eqref{EQ:THDonJohn2} follows.  \epr

Theorem~\ref{TH:DonJohn} shows that the normalizing function (rate)
$\psi_{n}^{01}(\theta)$ and the result~\eqref{low:3} cannot be
considered as a benchmark. Indeed, the risk of the {\sc es}
estimator is strictly below this bound. It attains the rate
$\psi_{n}^{01}(\theta)$ everywhere on $\R^n$
(cf.~\eqref{EQ:THDonJohn1}) and has strictly better rate on some
subsets of $\R^n$ (cf.~\eqref{EQ:THDonJohn2},~\eqref{EQ:infty}). In
particular, the {\sc es} estimator improves upon the soft
thresholding estimator, which is known to asymptotically attain the
bound~\eqref{low:3} \citep[cf.][]{DonJoh94a}. This is a kind of
inadmissibility statement for the rate $\psi_{n}^{01}(\theta)$.

Observe also that the improvement that we obtain is not a "marginal"
effect regarding signals $\theta$ with small intensity.
Indeed,~\eqref{EQ:infty} is stronger than~\eqref{EQ:THDonJohn2} and
the set $\Theta_n$ is rather massive. In particular, the $\ell_0$
norm $M(\theta)$ in the definition of $\Theta_n$ can be arbitrary,
so that $\Theta_n$ contains elements $\theta$ with the whole
spectrum of $\ell_1$ norms, from small $|\theta|_1=a/\sqrt{n}$ to
very large $|\theta|_1=bM/\sqrt{n}=b\sqrt{n}$. Various other
examples of $\Theta_n$ satisfying~\eqref{EQ:infty} can be readily
constructed.

So far, we were interested only in the rates. The fact that the
constant in~\eqref{EQ:THDonJohn1} is equal to~2 was of no importance
in this argument since on some subsets of $\R^n$ we can improve the
rate. Notice that one can construct estimators having the same
properties as those proved for $\tilde \theta^{\textsc{es}}$ in
Theorem~\ref{TH:DonJohn} with constant~1 instead of~2
in~\eqref{EQ:THDonJohn1}. In other words, one can construct an
estimator $\tilde \theta^*$ whose risk is at least as small as
$\psi_{n}^{01}(\theta)(1+o(1))$ everywhere on $\R^n$ and attains
strictly faster rate $o(\psi_{n}^{01}(\theta))$ on some subsets of
$\R^n$. Such an estimator $\tilde \theta^*$ can be obtained by
aggregating $\tilde \theta^{\textsc{es}}$ with the soft thresholding
estimator, as shown in the next theorem.

\begin{TH1}
\label{TH:DonJohn1} Consider the diagonal model.
 Then there exists a randomized estimator $\tilde\theta^*$
 such that
\begin{eqnarray}\label{EQ:THDonJohn11}
&&\limsup_{n\to\infty}\sup_{\theta\in \R^n} \frac{E_{\theta}|\tilde
\theta^*-\theta|_2^2 }{\tilde\psi_{n}(\theta)}\le 1,
\end{eqnarray}
where the expectation includes that over the randomizing
distribution, and where the normalizing functions $\tilde\psi_{n}$
satisfy
\begin{eqnarray}\label{EQ:THDonJohn12}
\liminf_{n\to\infty}\inf_{\theta\in \R^n}
\frac{\psi_{n}^{01}(\theta)}{\tilde\psi_{n}(\theta)}\ge
1, 
\end{eqnarray}
and
\begin{eqnarray}\label{EQ:THDonJohn13}
&&\liminf_{n\to\infty}\sup_{\theta\in \R^n}
\frac{\psi_{n}^{01}(\theta)}{\tilde\psi_{n}(\theta)} = \infty.
\end{eqnarray}
\end{TH1}
The proof of this theorem is given in the appendix.

\subsection{Minimax optimality on the intersection of $\ell_0$ and $\ell_1$ balls}
\label{SEC:opt}

The rate in the upper bound of Theorem~\ref{TH:sparseUB1} is the
minimum of terms depending on the $\ell_0$ norm $M(\theta)$ and on
the $\ell_1$ norm $|\theta|_1$, cf.~\eqref{EQ:psi}. We would like to
derive a corresponding lower bound, i.e., to show that this rate of
convergence cannot be improved in a minimax sense. Since both
$\ell_0$ and $\ell_1$ norms are present in the upper bound, a
natural approach is to consider minimax lower bounds on the
intersection of $\ell_0$ and $\ell_1$ balls. Here we prove such a
lower bound under some assumptions on the dictionary $\cH=\{f_1,
\ldots, f_M\}$ or, equivalently, on the matrix $\bX$. Along with the
lower bound for one ``worst case" dictionary $\cH$, we also state it
uniformly for all dictionaries in a certain class.

\subsubsection{Assumptions on the dictionary}\label{sec:ass_dict}

Recall first, that all the results from Section~\ref{SEC:UB} hold
under the only condition that the dictionary $\cH$ is composed of
functions $f_j$ such that $\|f_j\| \le 1$. This condition is very
mild compared to the assumptions that typically appear in the
literature on sparse recovery using $\ell_1$ penalization such as
the Lasso or the Dantzig selector. \citet{BuhGee09} review a long
list of such assumptions, including the restricted isometry (RI)
property given, for example, in \citet{Can08} and the restricted
eigenvalue (RE) condition of \citet{BicRitTsy09} described in
Section~\ref{SEC:lasso_bic}. We call them for brevity the
$L$-conditions. Loosely speaking, they ensure that for some integer
$S \le M$, the design matrix $\bX$ forms a quasi-isometry from a
suitable subset $\cA_{\pp}$ of $\R^\pp$ into $\R^n$ for any $\pp$
such that $|\pp|\le S$. Here ``quasi-isometry" means that there
exist two positive  constants $\underline \kappa$ and $\bar \kappa$
such that
\begin{equation}
\label{EQ:RI1}
\underline \kappa |\theta|_2^2 \le \frac{|\bX\theta|_2^2}{n} \le \bar \kappa |\theta|_2^2 \,,\qquad \forall\ \theta \in \cA_\pp\,.
\end{equation}
While the general thinking is that a design matrix $\bX$ satisfying
an $L$-condition is favorable, we establish below that, somewhat
surprisingly, such matrices correspond to the least favorable case.

We now formulate a weak version of the RI condition.
For any integer $M\ge 2$ and any $0<u\le M$ let
$\cP_u$ denote the set of vectors $\theta \in\{-1, 0, 1\}^M$ such
that $M(\theta)\le u$. For any constants
$\kappa\ge 1$  and $0<t\le (M\wedge n)/2$ let $\cD( t,
{ \kappa})$ 
be the class of design
matrices $\bX$ defined by the conditions:
\begin{itemize}
\item[(i)] $\DS \max_{1\le
j\le M} \|f_j\| \le 1$,
\item[(ii)] there exist $\underline{
\kappa}, \bar \kappa>0,$ such that $\underline{ \kappa}/ \bar \kappa\ge \kappa$ and
\begin{equation}
\label{EQ:class_dict}
 \underline{ \kappa}|\theta|_2^2\le  \frac{|\bX\theta|_2^2}{n}
\le \bar \kappa |\theta|_2^2, \qquad \forall\, \theta \in
\cP_{2t}\,.
\end{equation}
\end{itemize}%
Note that $t\le t'$ implies $\cD(t',\kappa)\le
\cD(t,\kappa) $. Examples of matrices $\bX$ that
satisfy~\eqref{EQ:class_dict} are given in the next subsection.

In the next subsection we show that the upper bound of
Theorem~\ref{TH:sparseUB1} matches a minimax lower bound which holds
uniformly over the class of design matrices~$\cD(S,\kappa)$. 

\subsubsection{Minimax lower bound} Denote by $P_\eta$ the distribution
of $(Y_1, \ldots, Y_n)$ where $Y_i=\eta(x_i)+\xi_i, \ i=1, \ldots,
n$, and by $E_\eta$ the corresponding expectation. For any
$\delta>0$ and any integers $S \ge 1, n\ge 1, M \ge 1, R\ge 1$ such that $R \le M \wedge n$, define the
quantity
\begin{equation} \label{EQ:zeta} \zeta_{n,M, R}( S,
\delta):=\frac{\sigma^2R}{n} \wedge \frac{\sigma^2S}{n}\log\left(1+\frac{eM}{S}\right)\wedge\frac{
\sigma\delta}{\sqrt n} \sqrt{\log \left(1+\frac{eM\sigma }{\delta
\sqrt n}\right)}\wedge \delta^2\,.
\end{equation}
Note that $\zeta_{n,M, R}( S, \delta)=\psi_{n,M}(\theta)$ where
$\psi_{n,M}$ is the function~\eqref{EQ:psi} with $
M(\theta)=S$ and $|\theta|_1=\delta$. Let $m \ge 1$ be the largest
integer satisfying
\begin{equation} \label{EQ:def_m}
m \le \frac{\delta \sqrt n}{\sigma\sqrt{\log\left(1+\frac{eM
}{m}\right)}}\,,
\end{equation}
if such an integer exists. If there is no $m\ge1$ such that~\eqref{EQ:def_m} holds, we set $m=0$. Note that $m \le
\delta\sqrt{n}/\sigma$.

\begin{TH1}\label{TH:low} Fix $\delta>0$ and integers $n\ge 1, M
\ge 2$, $1\le S\le M$. Fix $\kappa>0$ and let $\cH$ be
any dictionary with design matrix $\bX \in \cD(S\wedge \bar m, {
\kappa})$, where $\bar m=m\vee1$ and $m$ is defined
in~\eqref{EQ:def_m}. Then, for any estimator $T_n$, possibly
depending on $\delta, S, n, M$ and $\cH$, there exists a numerical
constant $c^*>0$, such that
\begin{equation}\label{low1}
\sup_{\substack{ \theta\in \R^M_+\setminus\{0\} \\ M(\theta)\le S \\
|\theta|_1\le\delta}}  \sup_{\eta} \left\{ E_\eta\|{T}_n -\eta\|^2 -
\| {\sf f}_{\theta}-\eta \|^2 \right\}\ge c^* \kappa
\zeta_{n,M, \rank{\bX}}(S, \delta)\,,
\end{equation}
where $\rank{\bX}$ denotes the rank of $\bX$ and $\R^M_+$ is the positive cone of $\R^M$. Moreover, 
\begin{equation}\label{low2}
\sup_{\substack{ \theta\in \R^M_+\setminus\{0\} \\
M(\theta)\le S \\ |\theta|_1\le\delta}} E_{\ff_\theta}\|T_n-
\ff_{\theta}\|^2\ge c^* \kappa \zeta_{n,M, \rank{\bX}}(S, \delta)\,.
\end{equation}
\end{TH1}
 The
proof of this theorem is given in Subsection~\ref{sub:proof_low} of
the appendix. It is worth mentioning that the result of Theorem~\ref{TH:low} is
stronger than the minimax lower bounds discussed in
Subsection~\ref{SEC:disc} (cf.~\eqref{low:2}) in the sense that even
if $\eta=\ff_{\theta^*}, \theta^*\in \R^M$, where $M(\theta^*)$ and
$|\theta^*|_1$ are known a priori, the rate cannot be improved.

%
%
%

Define $\widetilde R=
1+\left[\frac{\D}{C_0}\log\left(1+\frac{eM}{\D}\right)\right]$ for
some constant $C_0>0$ to be chosen small enough. We now show that
for each choice of $R\ge 1$ such that $\widetilde R \le M \wedge n$,
there exists at least one matrix $\bX \in \cD(R/2,
\kappa)$ such that $R \le \rank{\bX} \le \widetilde R$. A basic
example  is the following. Take the elements $\bX_{i,j}=f_j(x_i),
i=1, \ldots, n, j=1, \ldots, M$, of matrix $\bX$ as
\begin{equation}
\label{EQ:radmat}
\bX_{i,j}=\left\{
\begin{array}{ll}
\mathfrak{e}_{i,j} \sqrt{\frac{n}{\widetilde R}}
& \textrm{if}\ i\le \widetilde R\,,\\
0 & \textrm{otherwise }\,,
\end{array}\right.
\end{equation}
where $\mathfrak{e}_{i,j}, 1\le i \le, 1\le j \le M$ are i.i.d.
Rademacher random variables, i.e., random variables taking values
$1$ and $-1$ with probability $1/2$. First, it is clear that then
$\|f_j\|\le 1, j=1, \ldots, M$. Next, condition (ii) in the
definition of $\cD(R/2,\kappa)$ follows from the
results on RI properties of Rademacher matrices. Many such results
have been derived and we focus only on that of~\citet{BarDavDeV08}
because of its simplicity. Indeed, Theorem~5.2
in~\citet{BarDavDeV08} ensures not only that for an
integer $S'\le M\wedge n$ there exist design matrices in
$\cD(S'/2,\kappa)$ but also that most of the design matrices $\bX$
with i.i.d. Rademacher entries $\mathfrak{e}_{i,j}$ are in
$\cD(S'/2, \kappa)$ for some $\kappa>0$ as long as there exists a
constant $C_0$ small enough such that the condition
\begin{equation}
\label{EQ:condRIP2_a} \frac{ S'}{M\wedge
n}\log\left(1+\frac{eM}{S'}\right)< C_0\,
\end{equation}
is satisfied.  Specifically, Theorem~5.2 in~\citet{BarDavDeV08}
ensures that if $\bX'$ is the $\widetilde R\times M$ matrix composed
of the first $\widetilde R$ rows of $\bX$ with elements as defined
in (\ref{EQ:radmat}), and
\begin{equation}
\label{EQ:condRIP2_b} \frac{ \D}{\widetilde
R}\log\left(1+\frac{eM}{\D}\right)\le C_0\,
\end{equation}
holds for small enough $C_0$, then
\begin{equation*}
\label{EQ:RI1_a} \underline \kappa \frac{n}{\widetilde R}
|\theta|_2^2 \le \frac{|\bX'\theta|_2^2}{\widetilde R} \le \bar
\kappa \frac{n}{\widetilde R}|\theta|_2^2 \,,\qquad \forall\ \theta:
\, M(\theta)\le \D ,
\end{equation*}
with probability close to 1 which in turn implies (ii) with $t =
R/2$. As a result, the above construction yields $\bX
\in \cD(R/2, \kappa)$ that has rank bracketed by $R$
and $\widetilde R$ since~\eqref{EQ:condRIP2_b} holds by
our definition of $\widetilde R$.

In what follows  $C_0$ is the constant in~\eqref{EQ:condRIP2_a}
small enough to ensure that Theorem~5.2 in~\citet{BarDavDeV08}
holds, and we assume w.l.o.g. that $C_0<1$.

Using the above remarks and Theorem~\ref{TH:low} we obtain the
following result.
\begin{TH1}\label{COR:low}  Fix $\delta>0$ and integers $n\ge 1, M
\ge 2,  1\le S\le M, \D\ge 1$. Moreover, assume
that $1+\frac{R}{C_0}\log(1+eM/R) \le M \wedge n$. Then there
exists a dictionary $\cH$ composed of functions $f_j$ with
$\max_{1\le j\le M}\|f_j\|\le 1$, $R\le \rank{\bX}\le
1+\frac{R}{C_0}\log(1+eM/R)$, and a constant $c_*>0$ such that
\begin{equation}\label{low1_a}
\inf_{T_n}\sup_{\substack{ \theta\in \R^M_+\setminus \{0\} \\
 M(\theta)\le S \\ |\theta|_1\le\delta}}  \sup_{\eta}\left\{
E_\eta\|{T}_n -\eta\|^2 - \| {\sf f}_{\theta}-\eta \|^2 \right\}\ge
c_* \zeta_{n,M, \rank{\bX}}(S, \delta)\,.
\end{equation}
where the infimum is taken over all estimators.
Moreover, 
\begin{equation}\label{low2_a}
\inf_{T_n}\sup_{\substack{ \theta\in \R^M_+\setminus \{0\} \\
M(\theta)\le S \\ |\theta|_1\le\delta}} E_{\ff_\theta}\|T_n-
\ff_{\theta}\|^2\ge c_* \zeta_{n,M, \rank{\bX}}(S, \delta)\,.
\end{equation}
\end{TH1}
{\sc Proof.} Let $\bX$ be a random matrix constructed as
in~\eqref{EQ:radmat} so that the rank of $\bX$ is bracketed by $R$
and $\widetilde R$ and $\bX \in \cD(R/2, \kappa)$. We
consider two cases. Assume first that $S \le R/2$ so
that $\bX \in \cD(R/2, \kappa) \subset \cD(S,
\kappa)\subseteq \cD(S\wedge \bar m, \kappa)$ and the result follows
trivially from Theorem~\ref{TH:low}. Next, if $S\ge
R/2$, observe that
$$
\rank{\bX} \le  \tilde R \le 1 +
\frac{\D}{C_0}\log\left(1+\frac{eM}{\D}\right) \le
\frac{2}{C_0}R\log\left(1+\frac{2eM}{R}\right)\,,
$$
(we used here that $C_0<1$), so that
$$
\rank{\bX} \wedge S\log\left(1+\frac{eM}{S}\right) \le \rank{\bX}
\le  \frac{2}{C_0}\left(\rank{\bX} \wedge
R\log\left(1+\frac{2eM}{R}\right)  \right)\,.
$$
It yields $\zeta_{n,M, \rank{\bX}}(S, \delta) \le C \zeta_{n,M,
\rank{\bX}}( R/2, \delta)$ and the result follows from
Theorem~\ref{TH:low}, which ensures that
$$
\inf_{T_n}\sup_{\substack{ \theta\in \R^M_+\setminus \{0\} \\
 M(\theta)\le S \\ |\theta|_1\le\delta}}  \sup_{\eta}\left\{
E_\eta\|{T}_n -\eta\|^2 - \| {\sf f}_{\theta}-\eta \|^2 \right\}\ge
c^* \kappa \zeta_{n,M, \rank{\bX}}(R/2, \delta) \ge
c_*\zeta_{n,M, \rank{\bX}}(S, \delta)\,.
$$
\epr

As a consequence of Theorem~\ref{TH:low} we get a lower bound on the
$\ell_0$ ball $B_0(S)=\{\theta: M(\theta)\le S\}$  by formally
setting $\delta=\infty$ in~\eqref{low1}:
\begin{equation}\label{low1_0}
\sup_{\eta}\sup_{\substack{ \theta\in \R^M \\  M(\theta)\le S}}
\left\{  E_\eta\|{T}_n -\eta\|^2 - \| {\sf f}_{\theta}-\eta \|^2
\right\}\ge c_* \kappa
\frac{\sigma^2}{n}\left[\rank{\bX} \wedge S\log\left(1+\frac{eM}{S}\right)\right]
\end{equation}
and the same type of bound derived from~\eqref{low2}. Analogous
considerations 
lead to the following lower bound on the $\ell_1$ ball
$B_1(\delta)=\{\theta: |\theta|_1\le \delta\}$  when
setting $S=M$:
\begin{equation}\label{low1_1}
\begin{split}
\sup_{\eta}\sup_{\substack{ \theta\in \R^M \\ |\theta|_1\le\delta}}&
\left\{  E_\eta\|{T}_n -\eta\|^2 - \| {\sf f}_{\theta}-\eta \|^2
\right\}\ge c_* \kappa
\left(\frac{\sigma^2\rank{\bX}}{n}\wedge
\frac{\sigma\delta}{\sqrt n} \sqrt{\log \left(1+\frac{eM\sigma
}{\delta \sqrt n}\right)}\wedge \delta^2\right)\,,
\end{split}
\end{equation}%
and to the same type of bound derived from~\eqref{low2}. 

Consider now the linear regression, i.e., assume that there exists
$\theta^*$ such that $\eta={\sf f}_{\theta^*}$.
Comparing~\eqref{EQ:soiCORparma} with~\eqref{low2} we find that for
$\delta\ge 1/\sqrt{n}$ the rate $\zeta_{n,M, \rank{\bX}}(S, \delta)$ is
the minimax rate of convergence on $B_0(S)\cap B_1(\delta)$ and that
the {\sc es} estimator is rate optimal. Moreover, it is rate optimal
separately on $B_0(S)$ and $B_1(\delta)$, and the minimax rates on
these sets are given by the right hand sides of~\eqref{low1_0}
and~\eqref{low1_1} respectively.

For the  {\it diagonal model} (cf. Subsection~\ref{SEC:disc}),
asymptotic lower bounds and exact asymptotics of the minimax risk on
$\ell_q$ balls were studied by \citet{DonJohHoc92} for $q=0$ and  by
\citet{DonJoh94b} for $0<q<\infty$. These results were further
refined by \cite{AbrBenDon06}. In the $\ell_0$ case,
\citet{DonJohHoc92} exhibit a minimax rate over $B_0(S)$ that is
asymptotically equivalent to
$$
2\sigma^2\frac{S}{n}\log\left(\frac{n}{S}\right) \quad  {\rm as}\  M=n \to
\infty\,.
$$
In the $\ell_1$ case,~\citet{DonJoh94b} prove that the
minimax rate over an $\ell_1$ ball with radius $\delta$ is
asymptotically equivalent to
$$
\frac{\delta \sigma}{\sqrt
n}\sqrt{2\log\left(\frac{\sigma\sqrt{n}}{\delta}\right)}\quad  {\rm as}\  M=n
\to \infty\,.
$$
In both cases, the above rates are equivalent, up to a numerical
constant, to the asymptotics of the right hand sides
of~\eqref{low1_0} and~\eqref{low1_1} under the diagonal model. We
note that the results of those papers are valid under some
restrictions on asymptotical behavior of $S$ (resp. $\delta$) as a
function of $n$.

Recently \citet{RasWaiYu09} extended the study of asymptotic lower
bounds on $\ell_q$ balls ($0\le q \le 1$) to the non-diagonal case
with $M\ne n$. Their results hold under some restrictions on the
joint asymptotic behavior of $n, M$ and $S$
(respectively, $\delta$). 
The minimax rates on the $\ell_0$ and $\ell_1$ balls obtained in
\citet[Theorem~3]{RasWaiYu09} are similar to~\eqref{low1_0}
and~\eqref{low1_1} but, because of the specific asymptotics, some
effects are wiped out there. For example, the $\ell_1$ rate in
\cite{RasWaiYu09} is $\delta\sqrt{(\log M)/n}$,
whereas~\eqref{low1_1} reveals an elbow effect that translates into
different rates for $\sigma \rank{\bX}\le \delta \sqrt{n}$.
Furthermore, the dependence on the rank of $\bX$ does not appear in
\cite{RasWaiYu09}, since under their assumptions $\rank{\bX}=n$.
Theorem~\ref{TH:low} above gives a stronger result since it is (i)
non-asymptotic, (ii) it explicitly depends on the rank $\rank{\bX}$
of the design matrix and (iii) it holds on the intersection of the
$\ell_0$ and $\ell_1$ balls. Moreover, Theorem~\ref{TH:sparseUB1}
shows that the $\ell_0-\ell_1$ lower bound is attained by one single
estimator: the Exponential Screening estimator. Alternatively,
\citet{RasWaiYu09} treat the two cases separately, providing two
lower bounds and two different estimators that attain them in some
specific asymptotics.


\section{Universal aggregation}
\setcounter{equation}{0} \label{SEC:univ}

Combining the elements of a dictionary $\cH=\{f_1, \ldots, f_M\}$ to
estimate a regression function $\eta$ originates from the problem of
aggregation introduced by~\citet{Nem00}. It can be generally
described as follows. Given $\Theta \subset \R^M$, the goal of
aggregation is to construct an estimator $\hat f_n$ that satisfies
an oracle inequality of the form
\begin{equation}
\label{EQ:OIuniv}
\E\|\hat f_n -\eta\|^2 \le \min_{\theta \in \Theta} \|\ff_\theta-\eta\|^2 +C\Delta_{n,M}(\Theta)\,, \quad C>0\,,
\end{equation}
with the smallest possible (in a minimax sense) remainder term
$\Delta_{n,M}(\Theta)$, in which case $\Delta_{n,M}(\Theta)$ is
called optimal rate of aggregation, cf. \cite{Tsy03}. \citet{Nem00}
identified three types of aggregation: (MS) for \emph{model
selection}, (C) for \emph{convex} and (L) for \emph{linear}.
\citet{BunTsyWeg07c} also considered another collection of
aggregation problems, denoted by (${\rm L}_D$) for \emph{subset
selection} and indexed by $D \in \{1, \ldots, M\}$.  To each of
these problems corresponds a given set $\Theta \subset \R^M$ and an
optimal remainder term $\Delta_{n,M}(\Theta)$. For (MS) aggregation,
$\Theta=\Theta_{\rm (MS)}=B_0(1)\cap B_1(1)=\{e_1, \ldots, e_M\}$,
where $e_j$ is the $j$-th vector of the canonical basis of $\R^M$.
For (C) aggregation, $\Theta=\Theta_{\rm (C)}$ is a convex compact
subset of the simplex $B_1(1)=\{\theta \in \R^M\,:\, |\theta|_1\le
1\}$. The main example of $\{\ff_\theta, \theta \in \Theta_{\rm
(C)}\}$ is the set of all convex combinations the $f_j$'s. For (L)
aggregation, $\Theta=\Theta_{\rm (L)}=\R^M=B_0(M)$, so that
$\{\ff_\theta, \theta \in \Theta_{\rm (L)}\}$ is the set of all
linear combinations the $f_j$'s. Given an integer $D \in \{1,
\ldots, M\}$, for (${\rm L}_D$) aggregation, $\Theta=\Theta_{{\rm
(L}_D{\rm)}}=B_0(D)=\{\theta \in \R^M\,: M(\theta)\le D\}$. For this
problem, $\{\ff_\theta, \theta \in \Theta_{{\rm (L}_D{\rm)}}\}$ is
the set of all linear combinations  of at most $D$ of the $f_j$'s.

Note that all these sets $\Theta$ are of the form
$B_0(S)\cap B_1(\delta)$ for specific values of $S$ and $\delta$.
This allows us to apply the previous theory.

\renewcommand{\arraystretch}{1.5}
\begin{table}[h]
\begin{center}
\begin{tabular}{r|c|c}
Problem & $\Theta$ & $\Delta_{n,M}(\Theta)$ 
\\
\hline \hline (MS) & $\Theta_{\rm (MS)}=B_0(1)\cap B_1(1)$ &
$\frac{\log M}{n}$ 
\\
(C) & $\Theta_{\rm (C)}=B_1(1)$ & $\frac{M}{n}\wedge
\sqrt{\frac{1}{n}\log\left(1+\frac{eM}{\sqrt{n}}\right)}$
\\
(L) & $\Theta_{\rm (L)}=B_0(M)$ & $\frac{M}{n}$
\\
(${\rm L}_D$) & $\Theta_{{\rm (L}_D{\rm)}}=B_0(D)$ &
$\frac{D}{n}\log\left(1+\frac{eM}{D}\right)$
 \end{tabular}
\end{center}
\caption{Sets of parameters $\Theta_{\rm (MS)}, \Theta_{\rm (C)},
\Theta_{\rm (L)}$ and $\Theta_{{\rm (L}_D{\rm )}}$ and corresponding
optimal rates of aggregation presented in~\citet{BunTsyWeg07c}. Note
that~\citet{BunTsyWeg07c} considered a slightly different definition
in the (C) case: $\Theta_{\rm (C)}=B_1(1)\cap\R^M_+$ leading to the
same rate.} \label{TAB:theta}
\end{table}

Table~\ref{TAB:theta} presents the four different choices for
$\Theta$ together with the optimal remainder terms given by
\citet{BunTsyWeg07c}. For (MS), (C) and (L) aggregation they
coincide with optimal rates of aggregation originally proved in
\citet{Tsy03} for the regression model with i.i.d. random design and
integral $L_2$ norm in the risk. A fifth type of aggregation called
the $D$-{\it convex} aggregation, which we denote by (${\rm C}_{D}$)
was studied by \citet{Lou07}. In this case, $\Theta=\Theta_{{\rm
(C}_{D}{\rm )}}$ is a convex compact subset of $B_1(1)\cap B_0(D)$,
so that $\{\ff_\theta, \theta \in \Theta_{{\rm (C}_{D}{\rm )}}\}$
can be, as a typical example, the set of convex combinations of at
most $D$ of the $f_j$'s. \citet{Lou07} proves minimax lower bounds
together with an upper bound that departs from the lower bound by
logarithmic terms. However, the results hold in the i.i.d. random
design setting and do not extend to our setup. 
While several papers use different
estimators for different aggregation problems~\citep[see][]{Tsy03,
Rig09}, one contribution of~\citet{BunTsyWeg07c} was to show that
the {\sc bic} estimator defined in Section~\ref{SEC:lasso_bic}
satisfies oracle inequalities of the form
\begin{equation}\label{EQ:agg1}
\E\|\ff_{\hat \theta^\textsc{bic}} -\eta\|^2 \le (1+a)\min_{\theta \in \Theta} \|\ff_\theta-\eta\|^2 +C\frac{1+a}{a^2}\Delta_{n,M}(\Theta)\,,
\end{equation}
simultaneously for all the sets $\Theta$ presented in
Table~\ref{TAB:theta}. Here $a$ and $C$ positive constants.
Moreover, for the Lasso estimator defined in~\eqref{EQ:lasso},
\citet{BunTsyWeg07c} show less precise inequalities under the
assumption the matrix $\bX^\top \bX$ is positive definite, where
$\bX$ is the design matrix defined in Section~\ref{SEC:notation}.
Note that these oracle inequalities are not sharp since the leading
constant is $1+a$ and not 1, whereas letting $a\to0$ results in
blowing up the remainder term. The following theorem shows that the
Exponential Screening estimator satisfies sharp oracle inequalities
(i.e., with leading constant 1) that hold simultaneously for the
five problems of aggregation.
\begin{TH1}
\label{TH:univ} Assume that $\max_{1\le j\le M}\|f_j\|\le 1$. Then
for any $M \ge 2, n\ge 1, D\le M$, and $\Theta \in \{\Theta_{\rm
(MS)}, \Theta_{\rm (C)}, \Theta_{\rm (L)}, \Theta_{{\rm (L}_D{\rm
)}}, \Theta_{{\rm (C}_{D}{\rm )}}\}$ the Exponential Screening
estimator satisfies the following oracle inequality
\begin{equation*}
\E\|\ff_{\tilde \theta^{\textsc{es}}}-\eta\|^2 \le\min_{\theta \in
\Theta}\|\ff_{ \theta}-\eta\|^2 + C\Delta_{n,M}^*(\Theta)\,,
\end{equation*}
where $C>0$ is a numerical constant and 
\begin{equation*}
 \Delta_{n,M}^*(\Theta) =
\left\{\begin{array}{ll} \frac{ \sigma^2 \D}{n}\wedge \frac{
\sigma^2\log M}{n}
& {\rm if}\ \Theta=\Theta_{\rm (MS)},
\\
\frac{ \sigma^2
\D}{n}\wedge\sqrt{\frac{\sigma^2}{n}\log\left(1+\frac{eM\sigma}{\sqrt{n}}\right)}
& {\rm if}\ \Theta=\Theta_{\rm (C)},
\\
\frac{ \sigma^2 \D}{n} & {\rm if}\ \Theta=\Theta_{\rm (L)},
\\
\frac{ \sigma^2 \D}{n} \wedge \frac{ \sigma^2
D}{n}\log\left(1+\frac{eM}{D}\right) & {\rm if}\
\Theta=\Theta_{{\rm
(L}_{D}{\rm )}},
\\
\frac{ \sigma^2 \D}{n}
\wedge\sqrt{\frac{\sigma^2}{n}\log\left(1+\frac{eM\sigma}{\sqrt{n}}\right)}
\wedge \frac{\sigma^2 D}{n} \log\left(1+\frac{eM}{D }\right)
& {\rm if}\ \Theta=\Theta_{{\rm (C}_{D}{\rm )}}\,.\\
\end{array}
\right.
\end{equation*}
\end{TH1}
The proof of Theorem~\ref{TH:univ} follows directly
from~\eqref{EQ:soiCOR} and~\eqref{EQ:phi}.

 We also observe that $\Delta_{n,M}^*(\Theta)$ is to
within a constant factor of $\Delta_{n,M}^*(\Theta)\wedge 1$ since
$\frac{ \sigma^2 \D}{n}\le \frac{ \sigma^2 (M\wedge n)}{n} \le
\sigma^2$.

Using Theorems~\ref{TH:low} and~\ref{COR:low} it is not hard to show
that the rates $\Delta_{n,M}^*(\Theta)\wedge 1$ for
$\Delta_{n,M}^*(\Theta)$ listed in Theorem~\ref{TH:univ} are {\it
optimal rates of aggregation} in the sense of \citet{Tsy03}.
Indeed, it means to prove that there exists a dictionary $\cH$
satisfying the assumptions of Theorem~\ref{COR:low}, and a constant
$c>0$ such that the following lower bound holds:
\begin{equation}\label{low1_aa}
\inf_{T_n}\sup_{\eta}  \left\{  E_\eta\|{T}_n -\eta\|^2 -
\min_{\theta\in \Theta}\| {\sf f}_{\theta}-\eta \|^2 \right\}\ge c
(\Delta_{n,M}^*(\Theta)\wedge 1)\,,
\end{equation}
where the infimum is taken over all estimators. An important
observation here is that the left hand side of (\ref{low1_aa}) is
greater than or equal to
\begin{equation}\label{low2_aa}
\inf_{T_n}\sup_{ \theta\in \Theta} E_{\ff_\theta}\|T_n-
\ff_{\theta}\|^2.
\end{equation}
It remains to note that a lower bound
for~\eqref{low2_aa} with the rate $\Delta_{n,M}^*(\Theta)\wedge 1$
follows directly from Theorem~\ref{COR:low} (cf. also~\eqref{low1_0}
and~\eqref{low1_1}) applied with the values $S$ and $\delta$
corresponding to the definition of $\Theta$.

Interestingly, the rates given in Theorem~\ref{TH:univ} are
different from those in Table~\ref{TAB:theta}, and also from those
for the regression model with i.i.d. random design established in
\citet{Tsy03} and \citet{Lou07}. Indeed, they depend on the rank
$\D$ of the regression matrix $\bX$, and the bounds are better when
the rank is smaller. This is quite natural since the distance
$\|\ff_{\tilde \theta^{\textsc{es}}}-\eta\|$ is the ``empirical
distance" depending on $\bX$. One can easily understand it from the
analogy with the behavior of the ordinary least squares estimator,
cf. Lemma~\ref{LEM:lse}. Alternatively, the distance used in
\citet{Tsy03} and \citet{Lou07} for the i.i.d. random design setting
is the $L_2(P_X)$-distance where $P_X$ is the marginal distribution
of $X_i$'s, and no effects related to the rank can occur. As
concerns Table~\ref{TAB:theta}, the optimality of the rates given
there is proved in \citet{BunTsyWeg07c} only for $M\le n$ and
$\bX^\top\bX/n$ equal to the identity matrix, in which case $\D=M$
and thus the effect of $\D$ is not visible.

\renewcommand{\arraystretch}{1}

\section{Implementation and numerical illustration}
\label{SEC:simul} \setcounter{equation}{0} In this section, we
propose an implementation of the {\sc es} estimator together with a
numerical experiment both on artificial and real data. We suppose
throughout that the sample is fixed, so that the least squares
estimators $\theta_\pp, \pp\in\cP,$ are fixed vectors.

\subsection{Implementation via Metropolis approximation}
\label{SUB:MH}
Recall that the {\sc es} estimator ${\tilde \theta^{\textsc{es}}}$
is the following mixture of least squares estimators:
\begin{equation}
\label{EQ:defES}
{\tilde \theta^{\textsc{es}}}:=\frac{ \DS \sum_{\pp \in \cP} {\hat
\theta_\pp}\exp\Big(-\frac{1}{4\sigma^2}\sum_{i=1}^n (Y_i-\ff_{\hat
\theta_\pp}(x_i))^2-\frac{|\pp|}{2}\Big)\pi_{\pp}}{\DS \sum_{\pp \in
\cP}\exp\Big(-\frac{1}{4\sigma^2}\sum_{i=1}^n (Y_i-\ff_{\hat
\theta_\pp}(x_i))^2-\frac{|\pp|}{2}\Big)\pi_{\pp}}\,,
\end{equation}
where $\cP:=\{0, 1\}^M$, $\pi$ is the prior~\eqref{EQ:sprior}, and
$\hat \theta_\pp$ is the least squares estimator on $\R^\pp$. 

Recall also that the prior $\pi$ defined in~\eqref{EQ:sprior} assigns weight $1/2$ to the ordinary least squares estimator $\hat \theta_{\bone}$, where $\bone=(1, \ldots, 1) \in \cP$. It is not hard to check from the proof of Theorem~\ref{TH:sparseUB} that it allows us to cap the rates by $\sigma^2 R/n$. While this upper bound has important theoretical consequences, in the examples that we consider in this section, we typically have $R=n$ so that the dependence of the rates in $R$ is inconsequential. As a result, in the rest of the Section, we consider the following, simpler prior
\begin{equation}
\label{EQ:spriortilde}
\tilde \pi_\pp:=\left\{\begin{array}{ll}
\frac{2}{ H}\left(\frac{|\pp|}{2eM}\right)^{|\pp|}\,,& {\rm if \ } |\pp| < \D, \\
0\,,& {\rm otherwise}\,.
\end{array}\right.
\end{equation}
{\it Exact} computation of $\ff_{\tilde
\theta^{\textsc{es}}}$ requires the computation of $2^{R-1}$ least
squares estimators. In many applications this number is
prohibitively large and we need to resort to a numerical
approximation. Notice that $\tilde \theta^{\textsc{es}}$ is obtained
as the expectation of the random variable $\hat \theta_{\sf P}$
where ${\sf P}$ is a random variable taking values in $\cP$ with
probability mass function $\nu$ given by
$$
\nu_\pp\propto \exp\Big(-\frac{1}{4\sigma^2}\sum_{i=1}^n (Y_i-\ff_{\hat \theta_\pp}(x_i))^2-\frac{|\pp|}{2}\Big) \tilde \pi_{\pp}\,, \quad \pp \in \cP\,.
$$
This Gibbs-type distribution can be expressed as the stationary
distribution of the Markov chain generated by the
Metropolis-Hastings (MH) algorithm \citep[see,
e.g.,][Section~7.3]{RobCas04}. We now describe the MH algorithm
employed here. Consider the $M$-hypercube graph $\mathcal{G}$ with
vertices given by $\cP$. For any $\pp \in \cP$, define the
instrumental distribution $q(\cdot| \pp)$ as the uniform
distribution on the neighbors of $\pp$ in $\mathcal{G}$ and notice
that since each vertex has the same number of neighbors, we have
$q(\pp| \qq)=q(\qq| \pp)$ for any $\pp, \qq \in \cP$. The MH
algorithm is defined in Figure~\ref{FIG:MH}. We use here the uniform
instrumental distribution for the sake of simplicity. Our
simulations show that it yields satisfactory results both in
performance and in the speed. Another choice of $q(\cdot|\cdot)$ can
potentially further accelerate the convergence of the MH algorithm.
\begin{figure}[h]\label{FIG:MH}
\bookbox{
Fix $\pp_0=0\in \R^M$. For any $t\ge 0$, given $\pp_t \in \cP$,
\begin{enumerate}
\item Generate a random variable $\qQ_{t}$ with distribution
$q(\cdot|\pp_t)$.
\item Generate a random variable
$$
\pP_{t+1}=\left\{\begin{array}{lll}
\qQ_t & \textrm{with probability} & r(\pp_t, \qQ_t)\\
\pp_t & \textrm{with probability} & 1-r(\pp_t, \qQ_t)
\end{array}\right.
$$
where
$$
r(\pp, \qq)=\min\left( \frac{\nu_{\qq}}{\nu_{\pp}},1\right)\,.
$$
\item Compute the least squares estimator $\hat
\theta_{\pP_{t+1}}$.
\end{enumerate}}\caption{The Metropolis-Hastings algorithm on the $M$-hypercube.}
\end{figure}

The following theorem ensures the ergodicity of the Markov chain
generated by the MH algorithm.

\begin{TH1}
For any function $\pp \mapsto \theta_\pp \in \R^M$, the  Markov chain $(\pP_t)_{t\ge 0}$ defined by the MH algorithm satisfies
$$
\lim_{T\to \infty}\frac{1}{T}\sum_{t=T_0+1}^{T_0+T}\theta_{\pP_t}
=\sum_{\pp \in \cP}\theta_\pp \nu_\pp\,, \quad \nu-\text{almost surely}\,,
$$
where $T_0\ge 0$ is an arbitrary integer.
\end{TH1}
{\sc Proof.} The chain is clearly $\nu$-irreducible, so the result
follows from \citet[Theorem~7.4, p. 274]{RobCas04}. \epr

In view of this result, we approximate ${\tilde
\theta^{\textsc{es}}}=\sum_{\pp \in \cP}\theta_\pp \nu_\pp$ by
$$
\Tilde{\Tilde{\theta}}^{\textsc{es}}_T=
\frac{1}{T}\sum_{t=T_0+1}^{T_0+T}\hat
\theta_{\pP_t}\,,
$$
which is close to ${\tilde \theta^{\textsc{es}}}$ for sufficiently
large $T$. One salient feature of the MH algorithm is that it
involves only the ratios $\nu_\qq/\nu_\pp$ where $\pp$ and $\qq$ are
two neighbors in $\mathcal{G}$. Since
$$
\frac{\nu_{\qq}}{\nu_\pp}=\exp\Big(\frac{1}{4\sigma^2}
\sum_{i=1}^n\left[ (Y_i-\ff_{\hat \theta_\pp}(x_i))^2-(Y_i-\ff_{\hat \theta_\qq}(x_i))^2\right]+\frac{|\pp|-|\qq|}{2}\Big)\frac{ \tilde \pi_{\qq}}{\tilde \pi_\pp}\,,
$$
the MH algorithm benefits from the choice~\eqref{EQ:spriortilde} of the
prior $\tilde \pi$ in terms of speed. Indeed, for this prior, we have
$$
\frac{\tilde \pi_{\qq}}{\tilde \pi_\pp}=\left(1+\frac{\omega}{|\pp|}\right)^{|\qq|}\left(
\frac{|\pp|}{2eM}\right)^{\omega}\,,
$$
and $\omega=|\qq|-|\pp| \in \{-1, 1\}$ when $\pp$ and $\qq$ are two
neighbors in $\mathcal{G}$. In this respect, the choice of the prior
$ \tilde \pi$ as in~\eqref{EQ:spriortilde} is better than the suggestions
in~\citet{LeuBar06} and~\citet{Gir08} who consider priors that
require the computation of the combinatoric quantity ${M \choose
|\pp|}$. Moreover, the choice~\eqref{EQ:spriortilde} yields slightly
better constants and improves the remainder terms in the oracle
inequalities of Section~\ref{SEC:UB}, as compared to what would be
obtained with those priors.

As a result, the MH algorithm in this case takes the form of a
stochastic greedy algorithm with averaging, which measures a
tradeoff between sparsity and prediction to decide whether to add or remove a
variable. In all subsequent examples, we use a pure MATLAB implementation of the {\sc es} estimator. While the benchmark estimators considered below employ a C based code optimized for speed, we observed that a safe implementation of the MH algorithm (three time more iterations than needed) exhibited an increase of computation time of at most a factor two.

\subsection{Numerical experiments}

\subsubsection{Sparse recovery}

While our results for the {\sc es} estimator hold under no
assumption on the dictionary, we first compare the behavior of our
algorithm in a well-known example where the
$L$-conditions on the dictionary are satisfied and therefore sparse
recovery by $\ell_1$-penalized techniques is theoretically
achievable.

Consider the model $\bY=\bX\theta^*+\sigma \xi$, where $\bX$ is an
$n\times M$ matrix with independent Rademacher or standard Gaussian
entries and $\xi \in \R^n$ is a vector of independent standard Gaussian
random variables and is independent of $\bX$. The vector $\theta^*$
is given by $\theta^*_j=\1(j \le S)$ for some fixed $S$ so that
$M(\theta^*)=S$. The variance is chosen as $\sigma^2=S/9$ following
the numerical experiments of \citet[Section~4]{CanTao07}. For
different values of $(n,M,S)$, we run the {\sc es} algorithm on 500
replications of the problem and compare our results with several other popular estimators in
the sparse recovery literature. We limit our choice to estimators that are readily implemented in
R or MATLAB. The considered estimators are:
\begin{enumerate}
\item The
Lasso estimator with regularization parameter $\sigma\sqrt{8(\log
M)/n}$ as indicated in \citet{BicRitTsy09},
\item The cross-validated Lasso
estimator (LassoCV) with regularization parameter obtained by ten-fold
cross-validation,
\item The Lasso-Gauss
estimator (Lasso-G) corresponding to the Lasso estimator computed in 1., and threshold value given by $\sigma\sqrt{(2\log M)/n}$,
\item The cross-validated Lasso-Gauss
estimator (LassoCV-G) corresponding to the Lasso estimator computed in 2., and threshold value given by $\sigma\sqrt{(2\log M)/n}$,
\item The {\sc mc+} estimator of~\citet{Zha10} with regularization parameter $\sigma\sqrt{(2\log M)/n}$,
\item The {\sc scad} estimator  of~\citet{FanLi01} with regularization parameter $\sigma\sqrt{(2\log M)/n}$.
\end{enumerate}
The Lasso-Gauss
estimators in 3. and 4. are obtained using the
following two-step procedure. In the first step, a Lasso estimator (Lasso or LassoCV) is computed and only coordinates larger than the
threshold $\sigma\sqrt{2(\log M)/n}$ are retained in a set $\cJ$. In
the second step, the Lasso-Gauss estimators are obtained by
constrained least squares under the constraint that coordinates
$\hat \beta_j \notin \cJ$ are equal to 0. Indeed, it is usually
observed that the Lasso estimator induces a strong bias by over-shrinking large coefficients and the Lasso-Gauss procedure is a
practically efficient remedy to this issue. By construction, the {\sc scad} and {\sc mc+} estimators should not suffer from such a shrinkage. The Lasso estimators are based on the {\tt l1-ls} package in MATLAB \citep{KohKimBoy08}. The {\sc mc+} and {\sc scad} estimators are implemented in the {\tt plus} package in R \citep{ZhaMel09}.

\begin{figure}[h]
\psfrag{Index}[][][.8][0]{Index}
\psfrag{Iterations}[][][.8][0]{Iteration}
\psfrag{Value}[][][.8][0]{Value}
\begin{center}
\includegraphics[width=0.5\textwidth]{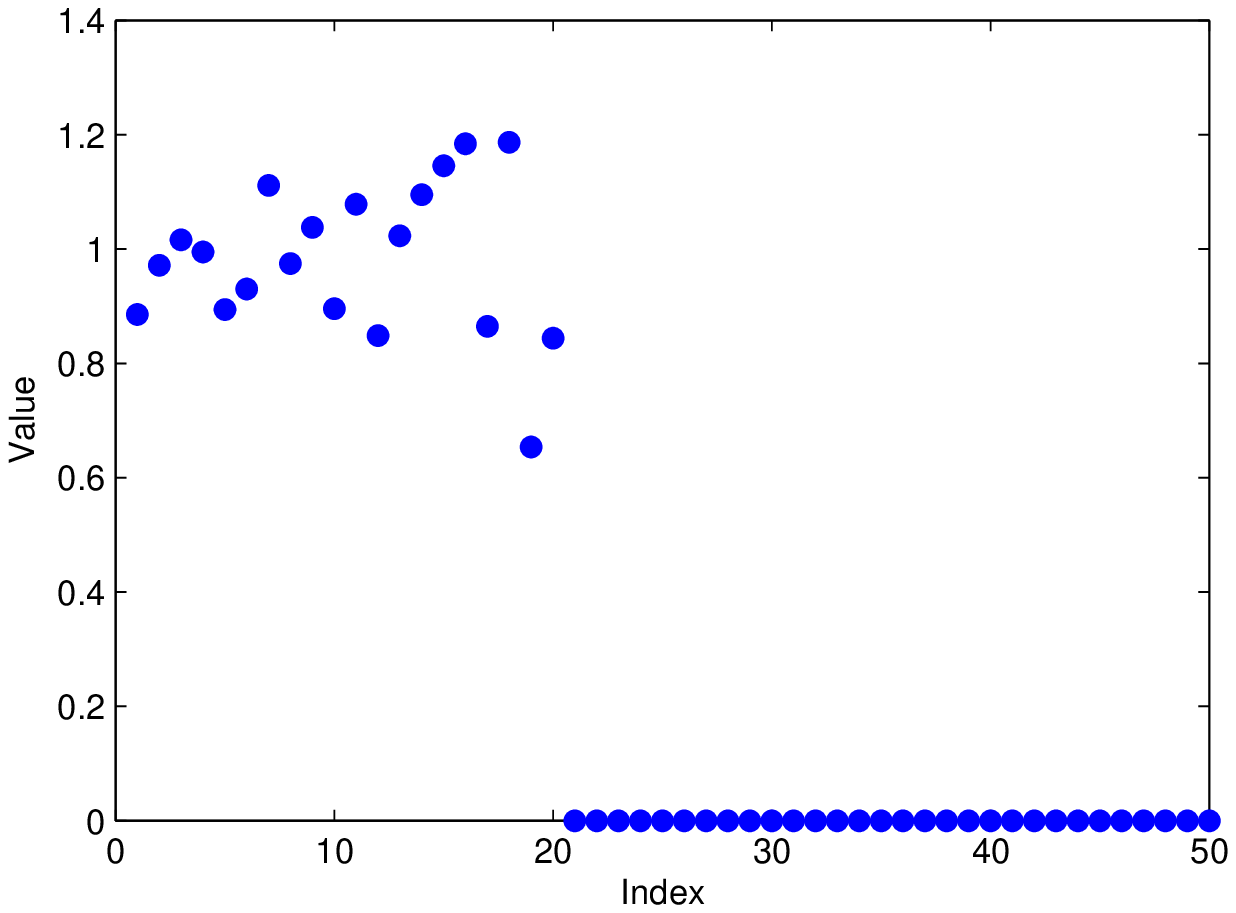}%
\includegraphics[width=0.5\textwidth]{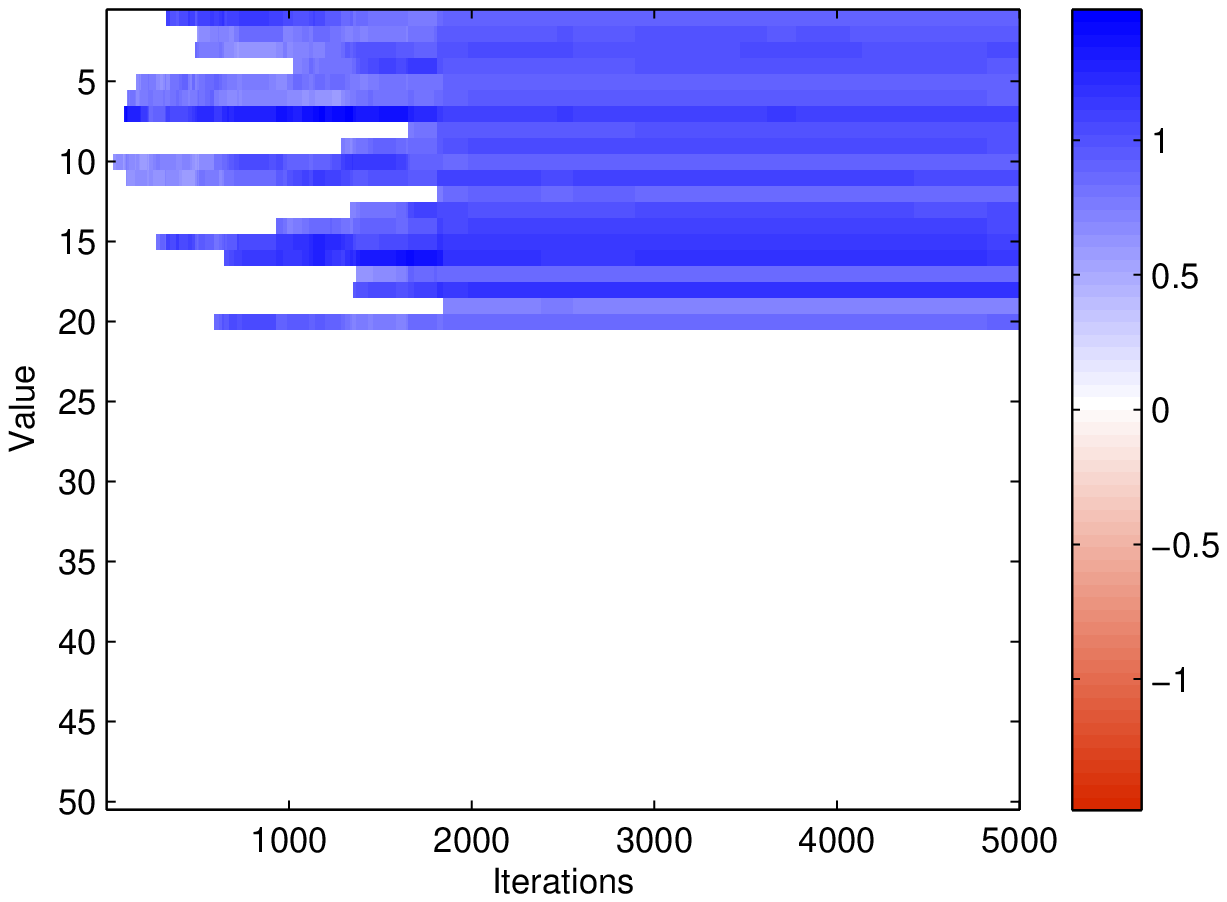}%
\end{center}
\caption{Typical realization for $(M,n,S)=(500,200,20)$ and Gaussian design. {\it Left:} Value of the  $\Tilde{\Tilde{\theta}}^{\textsc{es}}_T$, $T=7,000$, $T_0=3,000$. {\it Right:} Value of $\hat \theta_{\pP_t}$ for $t=1, \ldots, 5000$. Only the first 50 coordinates are shown for each vector.} \label{FIG:iter}
\end{figure}

The performance of each of the seven estimators generically denoted
by $\hat \theta$ is measured by its prediction error $|\bX(\hat
\theta -\theta^*)|_2^2/n=\|{\sf f}_{\hat \theta} -{\sf
f}_{\theta^*}\|^2$.
 Moreover,
even though the estimation error $|\hat \theta -\theta^*|_2^2$ is
not studied above, we also report its values in
Table~\ref{TAB:tabest}, for a better comparison with other
simulation studies. We considered the cases $(n,M,S)\in \{(100, 200,
10), (200, 500, 20)\}$. The Metropolis approximation
$\Tilde{\Tilde{\theta}}^{\textsc{es}}_T$ was computed with
$T_0=3,000$, $T=7,000$, which should be in the asymptotic regime of the Markov chain
since Figure~\ref{FIG:iter} shows that
on a typical example, the right sparsity pattern is recovered after
about 2,000 iterations.

\begin{figure}[h]
\psfrag{a}[tr][][.8][45]{{\sc es}}
\psfrag{b}[tr][][.8][45]{Lasso}
\psfrag{c}[tr][][.8][45]{LassoCV}
\psfrag{d}[tr][][.8][45]{Lasso-G}
\psfrag{e}[tr][][.8][45]{LassoCV-G}
\psfrag{f}[tr][][.8][45]{{\sc mc+}}
\psfrag{g}[tr][][.8][45]{{\sc scad}}
\begin{center}
\includegraphics[width=0.5\textwidth]{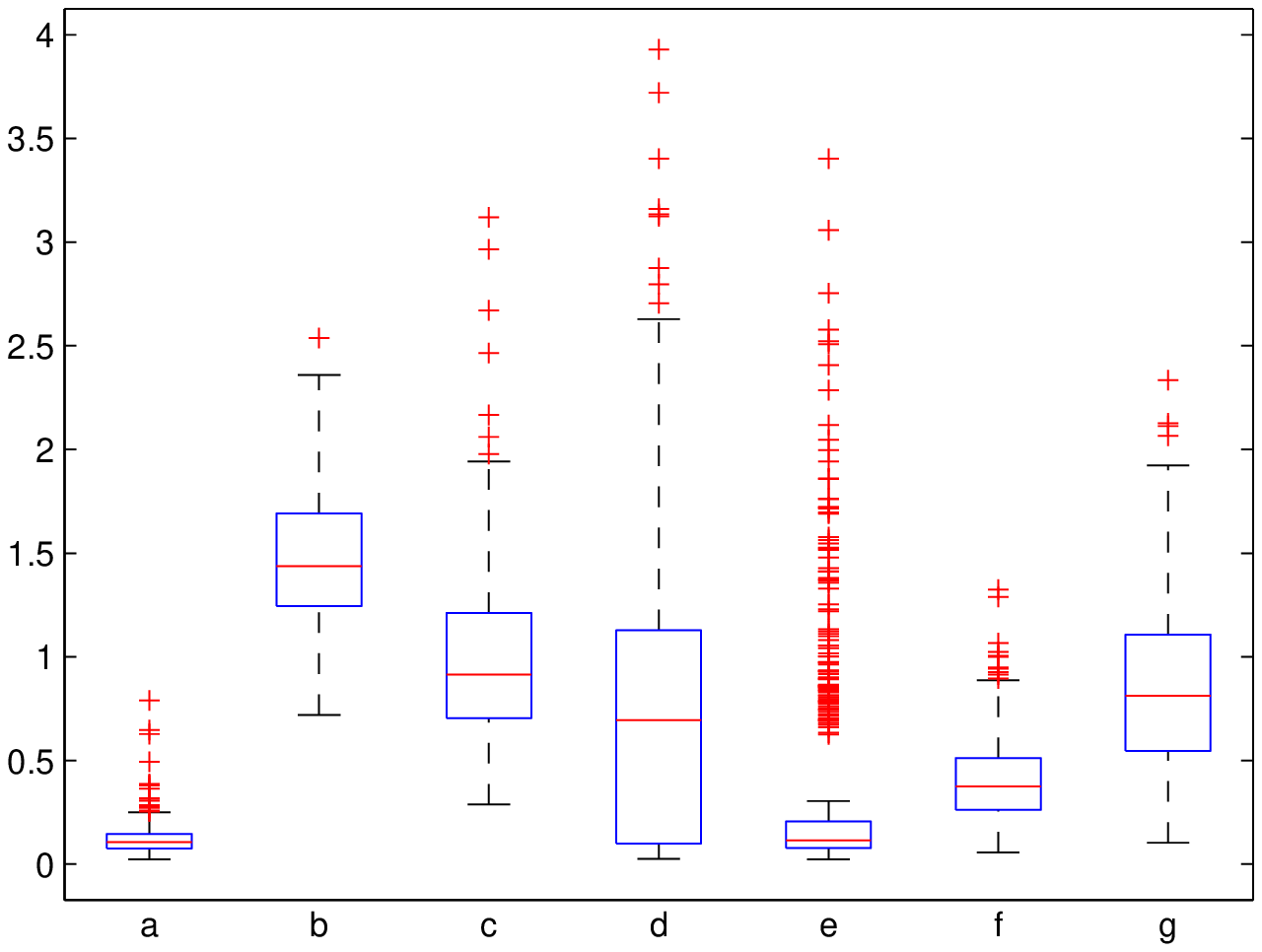}%
\includegraphics[width=0.5\textwidth]{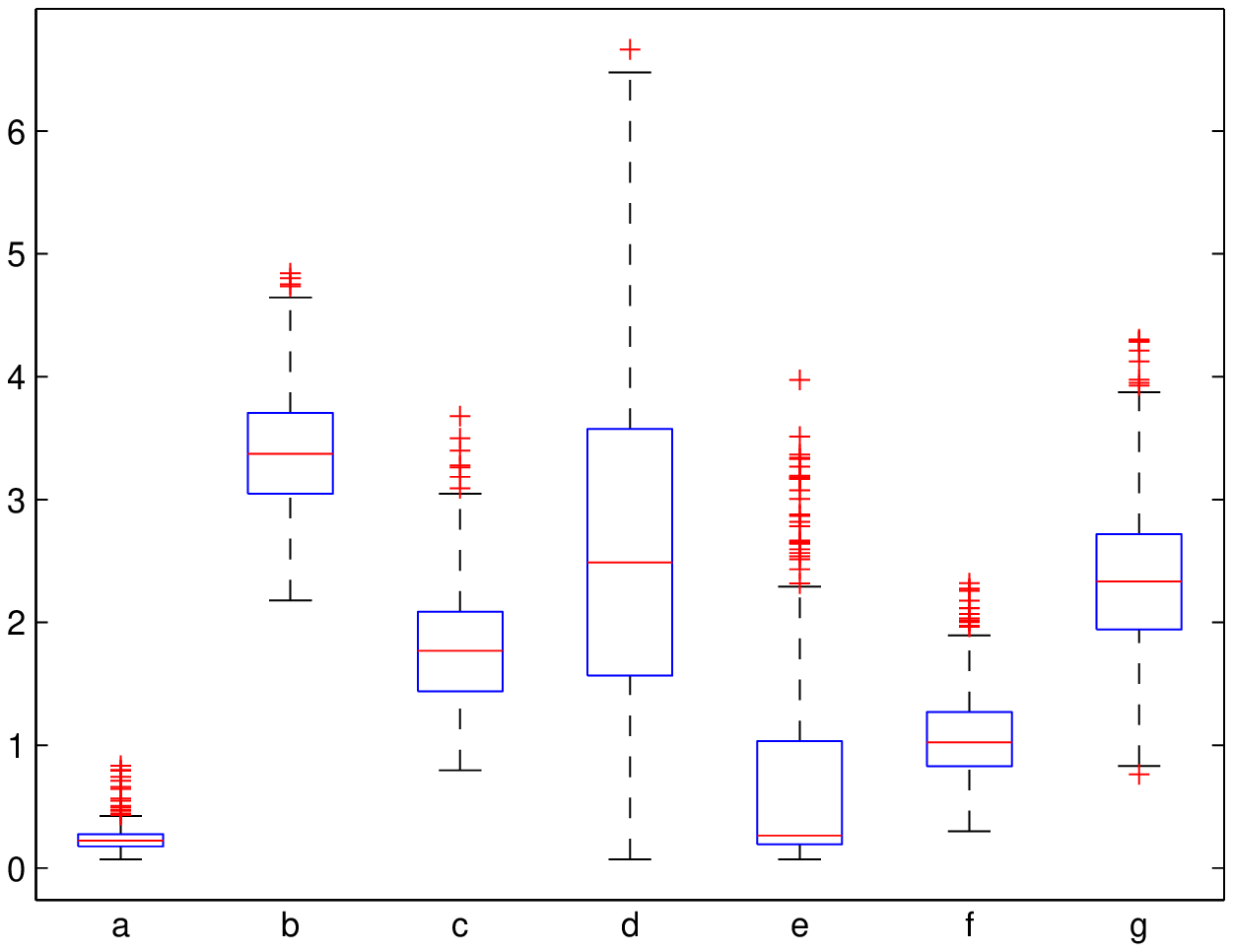}%
\end{center}
\vspace{-1.75cm}
\begin{center}
\includegraphics[width=0.5\textwidth]{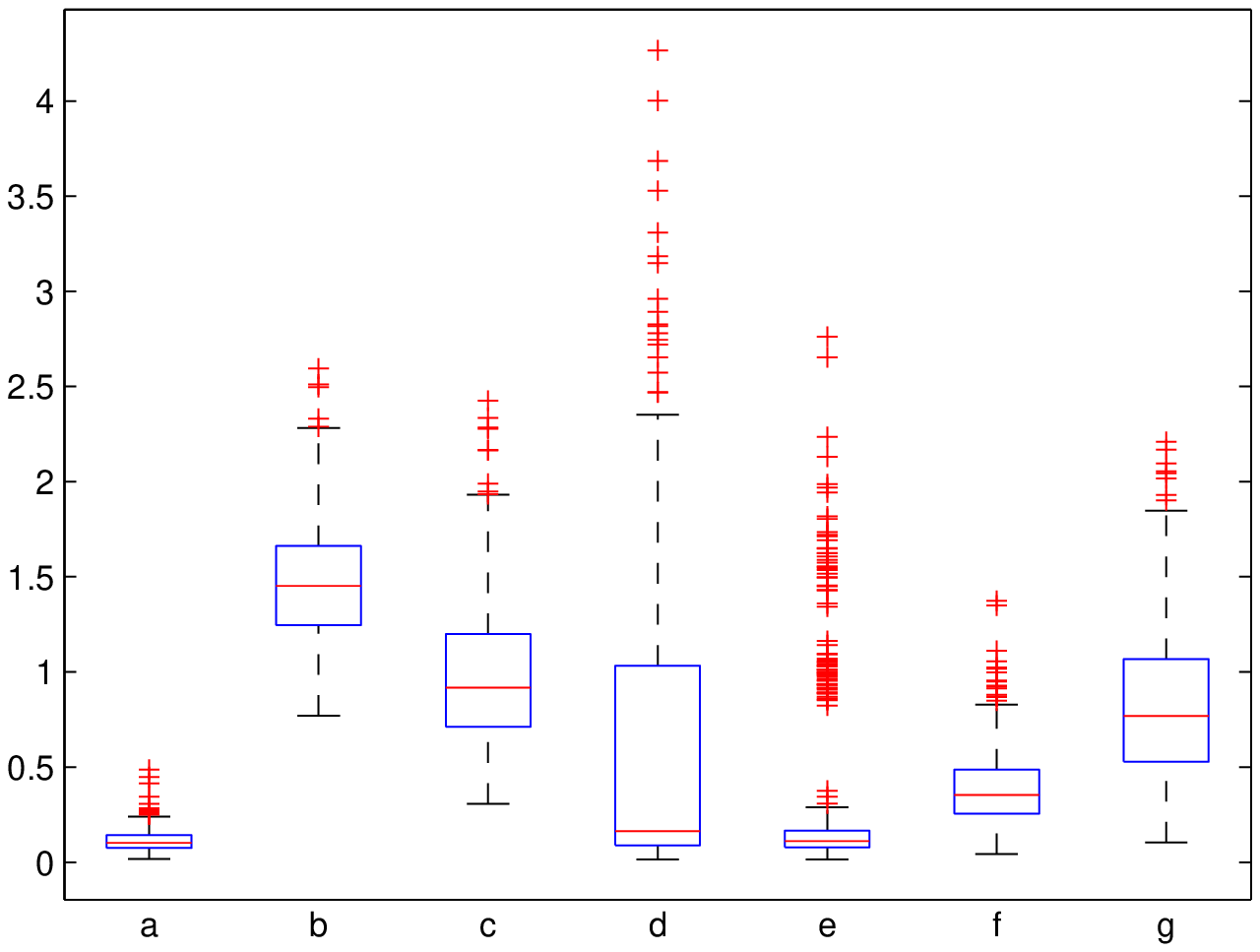}%
\includegraphics[width=0.5\textwidth]{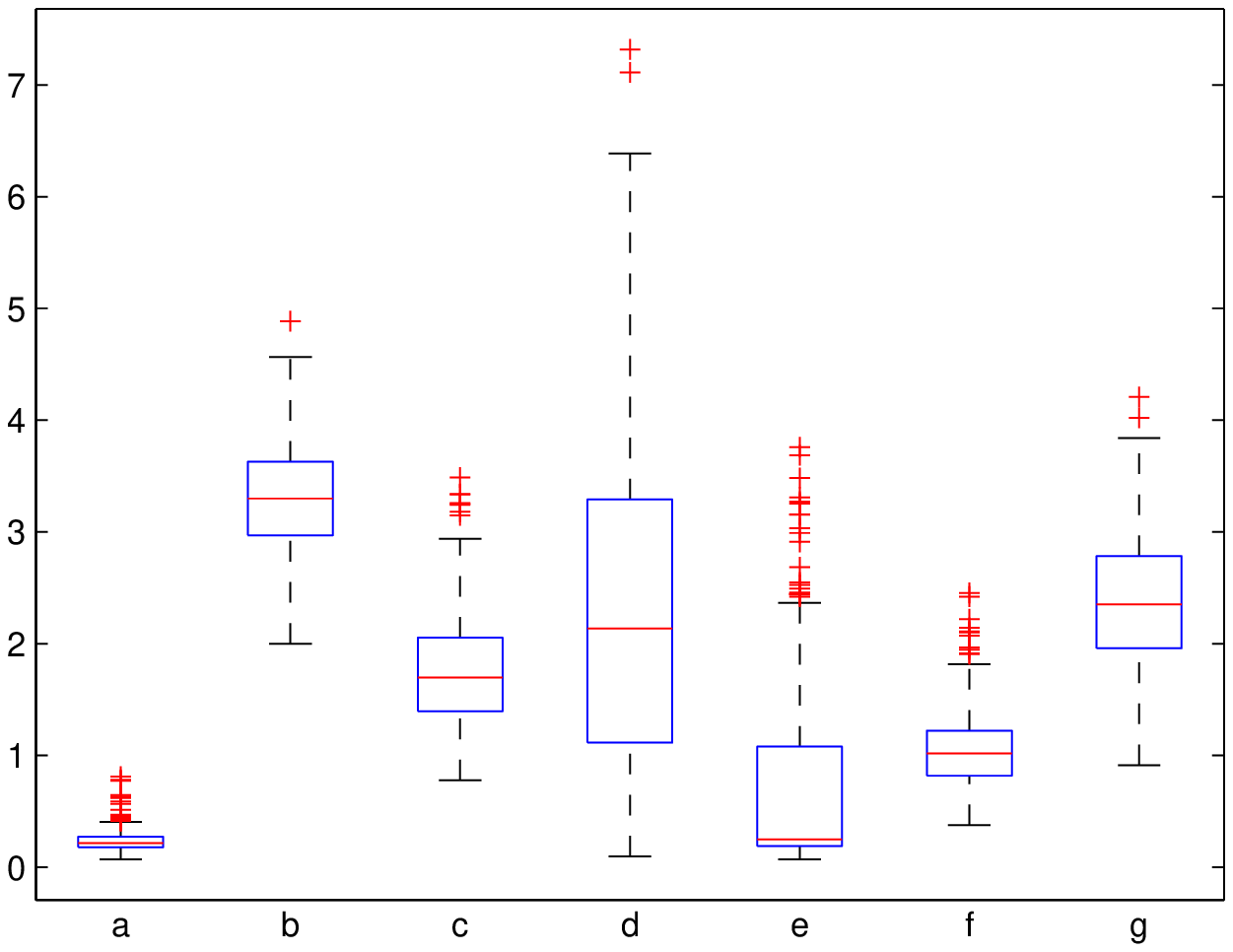}%
\end{center}
\caption{Boxplots of $|\bX(\hat \theta -\theta^*)|_2^2/n$ over 500 realizations for the {\sc es}, Lasso, cross-validated Lasso (LassoCV), Lasso-Gauss (Lasso-G), cross-validated Lasso-Gauss (LassoCV-G), {\sc mc+} and {\sc scad} estimators. {\it Left:} $(n,M,S)=(100,200,10)$, {\it right:} $(n,M,S)=(200,500,20)$, {\it top:} Gaussian design, {\it bottom:} Rademacher design.} \label{FIG:boxplots}
\end{figure}

Figure~\ref{FIG:boxplots} displays comparative boxplots for both
Gaussian and Rademacher design matrix. In particular, it shows that
{\sc es} outperforms all six other estimators and has less
variability across repetitions.

\begin{table}[h]
\begin{center}
\begin{tabular}{c|l|l|l|l|l|l|l}
$(M,n,S)$ & {\sc es} & Lasso & LassoCV & Lasso-G & LassoCV-G & {\sc mc+} & {\sc scad}\\
\hline
\hline
$(100,200, 10)$ & ${\bf 0.12}$  &  1.47  &  0.99   & 0.75 &   0.35 & 0.41 & 0.86  \\
& {\footnotesize (0.07)} &    {\footnotesize (0.31)}  &  {\footnotesize (0.40)  } & {\footnotesize (0.77) } &   {\footnotesize (0.53)} & {\footnotesize (0.20)} & {\footnotesize (0.40)}\\
\hline
$(200,500, 20)$ & {\bf  0.24}  &  3.39  &  1.81   & 2.55 &  0.70 & 1.07 & 2.37
\\
& {\footnotesize (0.10)} &    {\footnotesize (0.50)} &    {\footnotesize (0.50)}  &  {\footnotesize (1.45)  } & {\footnotesize (0.76) } &   {\footnotesize (0.35)}& {\footnotesize (0.64)} \\
\end{tabular}
\end{center}
\begin{center}
\begin{tabular}{c|l|l|l|l|l|l|l}
$(M,n,S)$ & {\sc es} & Lasso & LassoCV & Lasso-G & LassoCV-G & {\sc mc+} & {\sc scad}\\
\hline
\hline
$(100,200, 10)$ & ${\bf 0.12}$  &  1.48  &  0.99   & 0.70 &   0.30 & 0.39 & 0.83 \\
& {\footnotesize (0.06)} &    {\footnotesize (0.31)}  &  {\footnotesize (0.38)  } & {\footnotesize (0.79) } &   {\footnotesize (0.47)} &    {\footnotesize (0.19)}  &  {\footnotesize (0.39)  } \\
\hline
$(200,500, 20)$ & {\bf 0.24}  &  3.32  &  1.76  &  2.34  &   0.66 & 1.05 & 2.37
\\
& {\footnotesize (0.09)} &    {\footnotesize (0.49)}  &  {\footnotesize (0.49)  } & {\footnotesize (1.44) } &   {\footnotesize (0.74)}&    {\footnotesize (0.33)}  &  {\footnotesize (0.61)  }\\
\end{tabular}
\end{center}
\caption{Means and standard deviations of $|\bX(\hat \theta
-\theta^*)|_2^2/n$ over 500 realizations for the {\sc es}, Lasso, cross-validated Lasso (LassoCV), Lasso-Gauss (Lasso-G), cross-validated Lasso-Gauss (LassoCV-G), {\sc mc+} and {\sc scad} estimators. {\it Top:}
Gaussian design, {\it bottom:} Rademacher design.}\label{TAB:tabpred}
\end{table}

\begin{table}[h]
\begin{center}
\begin{tabular}{c|l|l|l|l|l|l|l}
$(M,n,S)$ & {\sc es} & Lasso & LassoCV & Lasso-G & LassoCV-G & {\sc mc+} & {\sc scad}\\
\hline
\hline
$(100,200, 10)$ & ${\bf 0.14}$  &  2.06  &  1.42   & 1.08 &   0.48 & 0.56 & 1.30  \\
& {\footnotesize {(0.12)}} &    {\footnotesize (0.72)}  &  {\footnotesize (0.66)  } & {\footnotesize (1.22) } &   {\footnotesize (0.84)} & {\footnotesize (0.34)} & {\footnotesize (0.81)}\\
\hline
$(200,500, 20)$ & {\bf  0.27}  &  4.72  &  2.73   & 3.62 &  0.93 & 1.45 & 3.51
\\
& {\footnotesize (0.13)} &    {\footnotesize (1.24)}  &  {\footnotesize (0.88)  } & {\footnotesize (2.29) } &   {\footnotesize (1.13)}& {\footnotesize (0.63)} &    {\footnotesize (1.33)} \\
\end{tabular}
\end{center}
\begin{center}
\begin{tabular}{c|l|l|l|l|l|l|l}
$(M,n,S)$ & {\sc es} & Lasso & LassoCV & Lasso-G & LassoCV-G & {\sc mc+} & {\sc scad}\\
\hline
\hline
$(100,200, 10)$ & ${\bf 0.13}$  &  1.99  &  1.37   & 0.94 &   0.38 & 0.51 & 1.21 \\
& {\footnotesize (0.07)} &    {\footnotesize (0.71)}  &  {\footnotesize (0.60)  } & {\footnotesize (1.19) } &   {\footnotesize (0.68)} &    {\footnotesize (0.35)}  &  {\footnotesize (0.81)  } \\
\hline
$(200,500, 20)$ & {\bf 0.26}  &  4.50  &  2.60  &  3.20  &   0.82 & 1.38 & 3.44
\\
& {\footnotesize (0.11)} &    {\footnotesize (1.14)}  &  {\footnotesize (0.80)  } & {\footnotesize (2.20) } &   {\footnotesize (1.00)}&    {\footnotesize (0.56)}  &  {\footnotesize (1.22)  }\\
\end{tabular}
\end{center}
\caption{Means and standard deviations of $|\hat \theta
-\theta^*|_2^2$ over 500 realizations for the {\sc es}, Lasso, cross-validated Lasso (LassoCV), Lasso-Gauss (Lasso-G), cross-validated Lasso-Gauss (LassoCV-G), {\sc mc+} and {\sc scad} estimators. {\it Top:}
Gaussian design, {\it bottom:} Rademacher design.}\label{TAB:tabest}
\end{table}

Figure~\ref{FIG:iter} illustrates a typical behavior of the
{\sc es} estimator for one particular realization of $\bX$ and
$\xi$. For better visibility, both displays represent only the 50
first coordinates of  $\Tilde{\Tilde{\theta}}^{\textsc{es}}_T$,
with $T=7,000$ and $T_0=3,000$. The left hand side display shows
that the sparsity pattern is well recovered and the estimated values
are close to one. The right hand side display illustrates the
evolution of the intermediate parameter $\hat \theta_{\pP_t}$ for
$t=1, \ldots, 5000$. It is clear that the Markov chain that runs on
the $M$-hypercube graph gets \emph{trapped} in the vertex that
corresponds to the sparsity pattern of $\theta^*$ after only $2,000$
iterations. As a result, while the {\sc es} estimator is not sparse itself, the MH approximation to the {\sc es} estimator may output a sparse solution.  A covariate $X_j$ is considered to be \emph{selected} by an estimator $\hat \theta$, if $|\hat \theta_j|>1/n$. Hence, for any two vectors $\theta^{(1)}, \theta^{(2)} \in \R^M$ define $\theta^{(1)}\symdiffsmall \theta^{(2)} \in \{0,1\}^M$ as the binary vector with $j$-th coordinate given by
$$
(\theta^{(1)} \symdiffsmall \theta^{(2)})_j=
\1(|\theta^{(1)}_j|>1/n, \theta^{(2)}_j=0)+ \1(\theta^{(1)}_j=0,
|\theta^{(2)}_j|>1/n)\,.
$$
The performance of an estimator $\hat \theta$ in terms of model selection is measured by the number $M(\hat \theta \symdiffsmall \theta^*)$ of variables that are incorrectly selected or incorrectly left out of the model. Among the four procedures considered here, {\sc mc+} uniformly dominates the other three in terms of model selection. Table~\ref{TAB:MS} displays the relative average model selection error (RAMS) over 500 repetitions of each of the experiments described above:
\begin{equation}
\label{EQ:rams}
\textrm{RAMS}(\hat \theta)=\frac{\sum_{i=1}^{500} M(\hat \theta^{(i)} \symdiffsmall \theta^*)}{\sum_{i=1}^{500} M(\hat \theta^{(i), \textsc{mc+}} \symdiffsmall \theta^*)}\,,
\end{equation}
where for each repetition $i$ of the experiment, $\hat \theta^{(i),\textsc{mc+}}$ denotes the \textsc{mc+} estimator and $\hat \theta^{(i)}$ is one of the four estimators: \textsc{es}, Lasso, \textsc{mc+} or \textsc{scad}.

\begin{table}[h]
\begin{center}
\begin{tabular}{c|l|l|l|l}
Design $(M,n,S)$ & {\sc es} & Lasso & {\sc mc+} & {\sc scad}\\
\hline
\hline
Gauss. $(100,200, 10)$ & 10.54  &  12.43  &  1.00   & 3.56  \\
\hline
Gauss. $(200,500, 20)$ & 9.26  &  15.81  &  1.00   & 6.04\\
\hline
\hline
Rad. $(100,200, 10)$ & 13.18    & 15.80  &  1.00  &  3.59 \\
\hline
Rad. $(200,500, 20)$ & 10.07  &  16.18  & 1.00    & 6.18
\\\end{tabular}
\end{center}
\caption{Relative average model selection error (RAMS) defined in~\eqref{EQ:rams} over 500 realizations for the {\sc es}, Lasso, {\sc mc+} and {\sc scad} estimators. {\it Top:}
Gaussian design, {\it bottom:} Rademacher design.}\label{TAB:MS}
\end{table}

While {\sc mc+} uniformly dominates the three other procedures, the
model selection properties of {\sc es} are better than Lasso but not
as good as {\sc scad} and the relative performance of {\sc es}
improves when the problem size increases. The superiority of {\sc
mc+} and {\sc scad} does not come as a surprise as these procedures
are designed for variable selection. However, {\sc es} makes up for
this deficiency by having much better estimation and prediction
properties.

\begin{figure}[h]
\psfrag{a}[tc][][.8][0]{Gauss. (100, 200, 10)}
\psfrag{b}[tc][][.8][0]{Rad. (100, 200, 10)}
\psfrag{c}[tc][][.8][0]{Gauss. (200, 500, 20)}
\psfrag{d}[tc][][.8][0]{Rad. (200, 500, 20)}
\begin{center}
\includegraphics[width=0.5\textwidth]{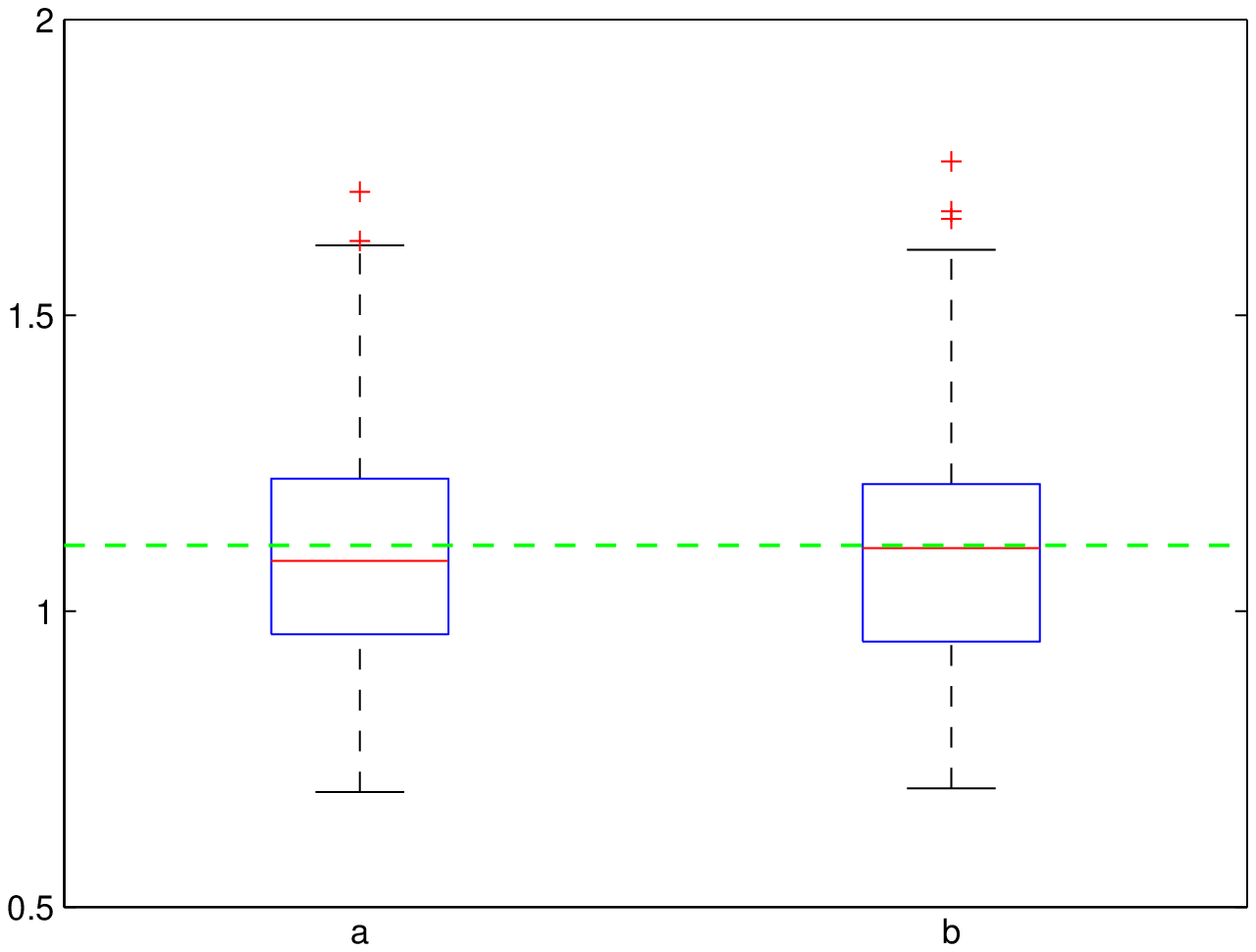}%
\includegraphics[width=0.5\textwidth]{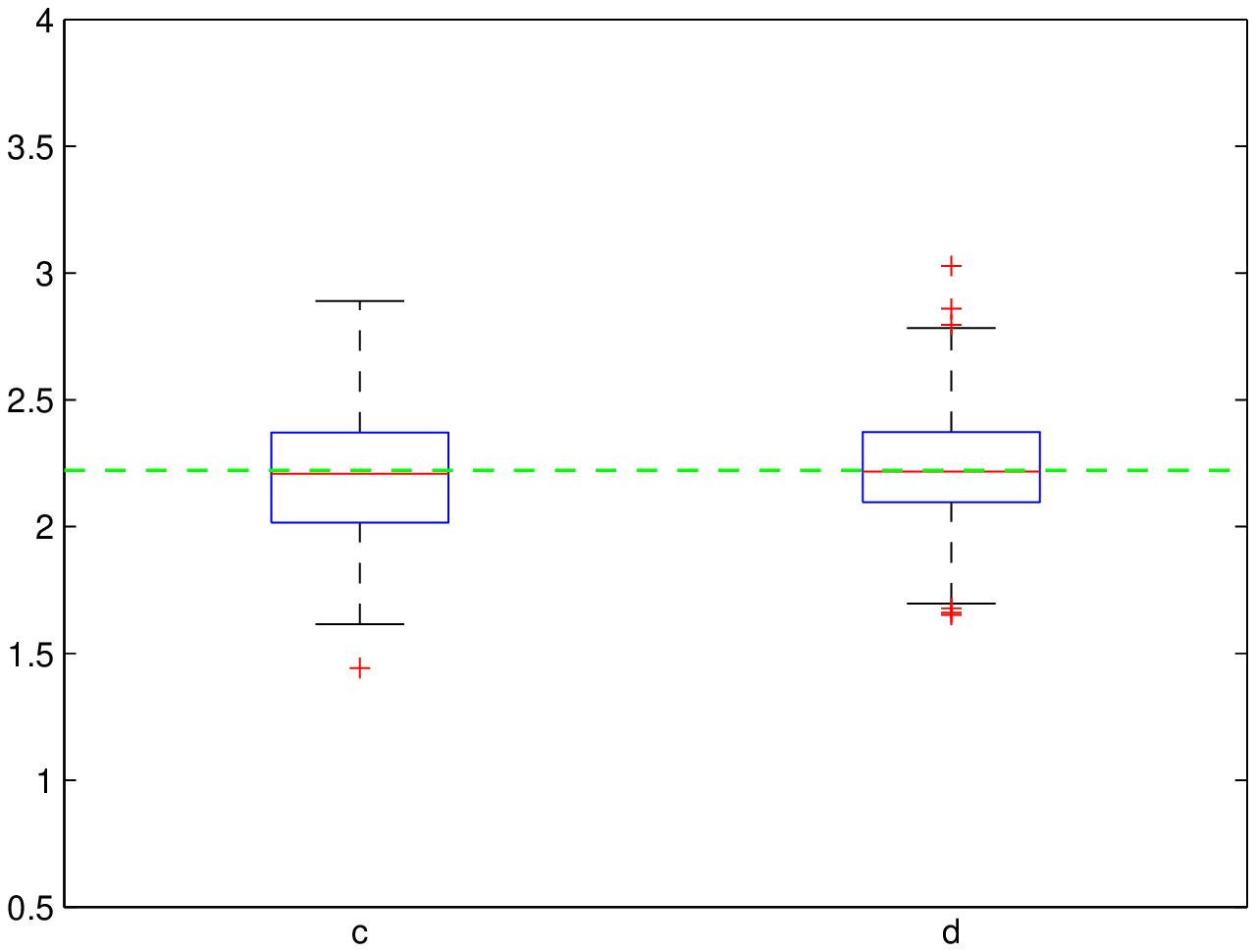}%
\end{center}
\caption{Boxplots of the estimated variance $\hat \sigma^2$ based on 500 replications of each of the four experiments described above. The horizontal dashed lines indicate the value of the true parameter $\sigma^2=S/9$. {\it Left:} $\sigma^2=1.11$.  {\it Right:} $\sigma^2=2.22$.} \label{FIG:varest}
\end{figure}

To conclude this numerical experiment in the linear regression model, notice that we used the knowledge of the variance parameter~$\sigma^2$ to construct the estimators, except for those based on cross-validation. In particular, \textsc{es} depends on $\sigma^2$ and it necessary to be able to implement it without such a knowledge. While an obvious solution consists in resorting to cross-validation or bootstrap, such procedures tend to become computationally burdensome. We propose the following estimator for $\sigma^2$. Let $\bar \theta^{\textsc{es}}$ denote the estimator obtained by replacing $\sigma^2$ with any upper bound $\bar \sigma^2 \ge \sigma^2$ in the definition~\eqref{EQ:defES} of the \textsc{es} estimator. Define
$$
\hat \sigma^2 =\inf \left\{ s^2: \left|\frac{|\bY - \bX\bar \theta^{\textsc{es}}(s^2)|^2_2}{n-M_n(\bar \theta^{\textsc{es}}(s^2))} - s^2 \right|> \alpha\right\}\,,
$$
where $\alpha>0$ is a tolerance parameter and for any $\theta \in
\R^M$, $M_n(\theta)=\sum_{j=1}^M \1(|\theta_j|>1/n)$. As a result,
the proposed estimator $\hat \sigma^2$ is the smallest positive
value that departs from the usual estimator for the variance by more
than $\alpha$. The motivation for this estimator comes from the
following heuristics, which is loosely inspired
by~\citet[Section~5.2]{Zha10}. It follows from the results
of~\citet{LeuBar06} that $\bar \theta^{\textsc{es}}(\bar \sigma^2)$
satisfies the oracle inequalities of Section~\ref{SEC:UB} and thus
of Section~\ref{SEC:univ} with $\sigma^2$ replaced by $\bar
\sigma^2$. As a consequence, we can use any  upper bound $\bar
\sigma^2 \ge\sigma^2$ to compute an estimator $\bar
\theta^{\textsc{es}}(\bar \sigma^2)$ and thus, an estimator of the
variance based on the residuals. Our heuristics consists in choosing
the smallest upper bound that is inconsistent with the estimator
based on the residuals. Figure~\ref{FIG:varest} and
Table~\ref{TAB:varest}  summarize the performance of the variance
estimator $\hat \sigma^2$ and the corresponding {\sc es} estimator
$\bar \theta^{\textsc{es}}(\hat \sigma^2)$ for $\alpha=1$.

\begin{table}[h]
\begin{center}
\begin{tabular}{c|l|l}
Design $(M,n,S)$ & $|\bX(\bar \theta^{\textsc{es}}(\hat \sigma^2)
-\theta^*)|_2^2/n$ & $|\bar \theta^{\textsc{es}}(\hat \sigma^2)
-\theta^*|_2^2$\\
\hline
\hline
Gauss. $(100,200, 10)$ & 0.12 & 0.14  \\
& {\footnotesize (0.09)} &    {\footnotesize (0.14)}  \\
\hline
Gauss. $(200,500, 20)$ & 0.26  &  0.31 \\
& {\footnotesize (0.19)} &    {\footnotesize (0.32)}  \\
\hline
\hline
Rad. $(100,200, 10)$ & 0.12 & 0.13 \\
& {\footnotesize (0.07)} &    {\footnotesize (0.08)}  \\
\hline
Rad. $(200,500, 20)$ & 0.25 & 0.28 \\
& {\footnotesize (0.11)} &    {\footnotesize (0.14)}  \\
\end{tabular}
\end{center}
\caption{Means and standard deviations of prediction error $|\bX(\bar \theta^{\textsc{es}}(\hat \sigma^2)
-\theta^*)|_2^2/n$ and estimation error $|\bar \theta^{\textsc{es}}(\hat \sigma^2)
-\theta^*|_2^2$  over 500 realizations for the {\sc es} estimator $\bar \theta^{\textsc{es}}(\hat \sigma^2)$ with estimated variance. \label{TAB:varest}}
\end{table}
Notice that in Table~\ref{TAB:varest}, the obtained values are
comparable to those in Tables~\ref{TAB:tabpred}
and~\ref{TAB:tabest}. It is worth noticing that the experiment with
Gaussian design and $(M,n,S)=(200,500,20)$ suffers from a long tail
of relatively poor performance (30 realizations out of 500 are
outliers) that deteriorates both the average performance and its
standard deviation. Nevertheless, it is remarkable that the {\sc es}
estimator with such estimator of the variance still has smaller
prediction and estimation errors in these experiments than the other
six considered methods.

\subsubsection{Handwritten digits dataset}

\begin{figure}[h]
\psfrag{cor}[tc][][.8][0]{{Correlation coefficient}}
\psfrag{count}[bc][][.8][0]{Count}
\begin{center}
\includegraphics[width=0.5\textwidth]{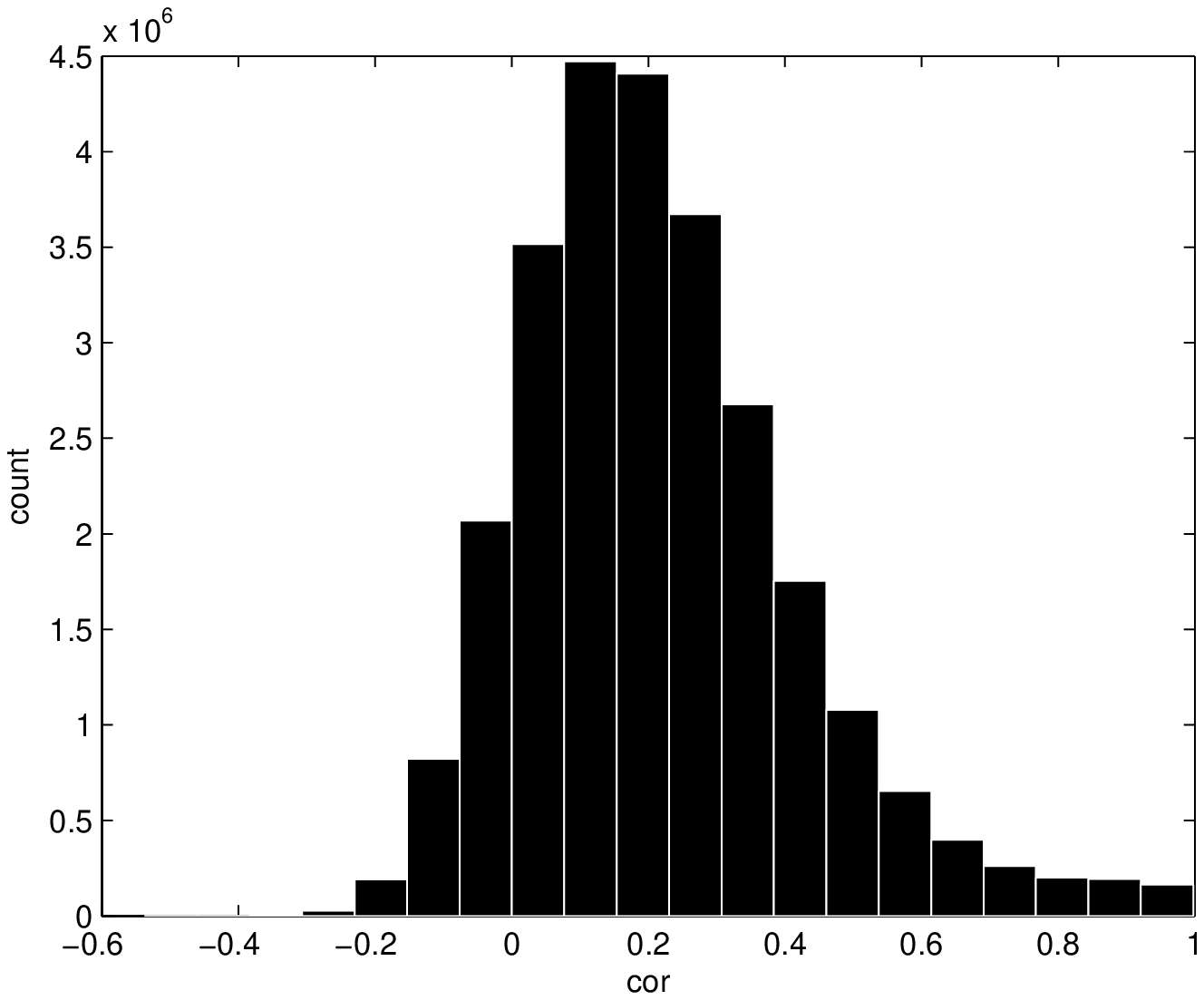}%
\includegraphics[width=0.5\textwidth]{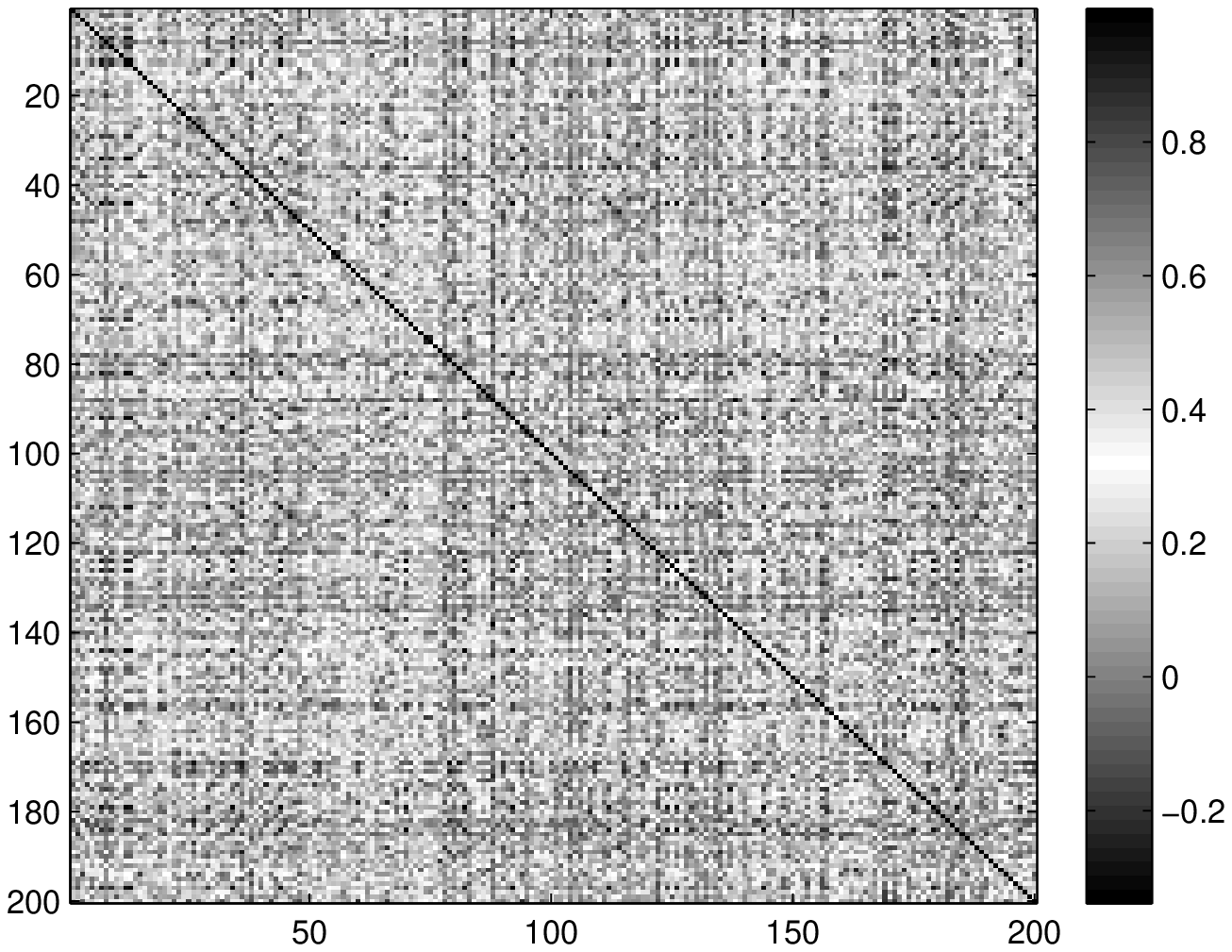}%
\end{center}
\caption{{\it Left:} Histogram of the $M(M-1)/2$ correlation coefficients between different images in the database.  {\it Right:} The upper left corner of size $200\times200$ of the full correlation matrix. Notice that only the absolute value of the correlation coefficients is discriminative in terms of color. The dark, off-diagonal regions are characteristic of correlated features.} \label{FIG:cor_digits}
\end{figure}

\begin{figure}[h] \centering
\subfigure[True]{%
\includegraphics[width=0.2\textwidth]{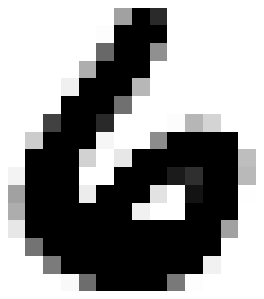}}%
\subfigure[Noisy]{%
\includegraphics[width=0.2\textwidth]{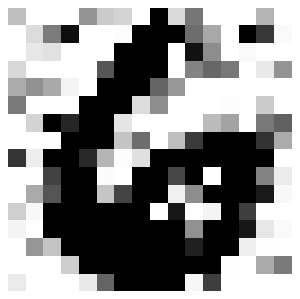}}%
\subfigure[{\sc es}]{%
\includegraphics[width=0.2\textwidth]{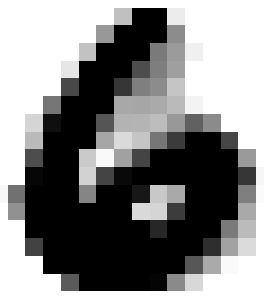}}%
\subfigure[Lasso]{%
\includegraphics[width=0.2\textwidth]{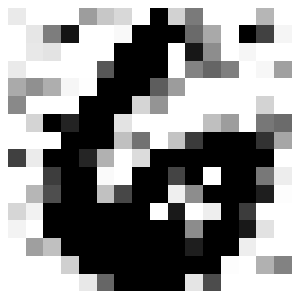}}%
\subfigure[Lasso-Gauss]{%
\includegraphics[width=0.2\textwidth]{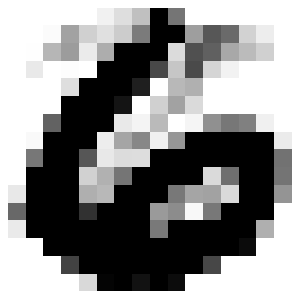}}%
\caption{Reconstruction of the digit ``6" with $\sigma=0.5$.} \label{FIG:6_05}
\end{figure}

\begin{figure}[t]
\begin{center}
\includegraphics[width=0.5\textwidth]{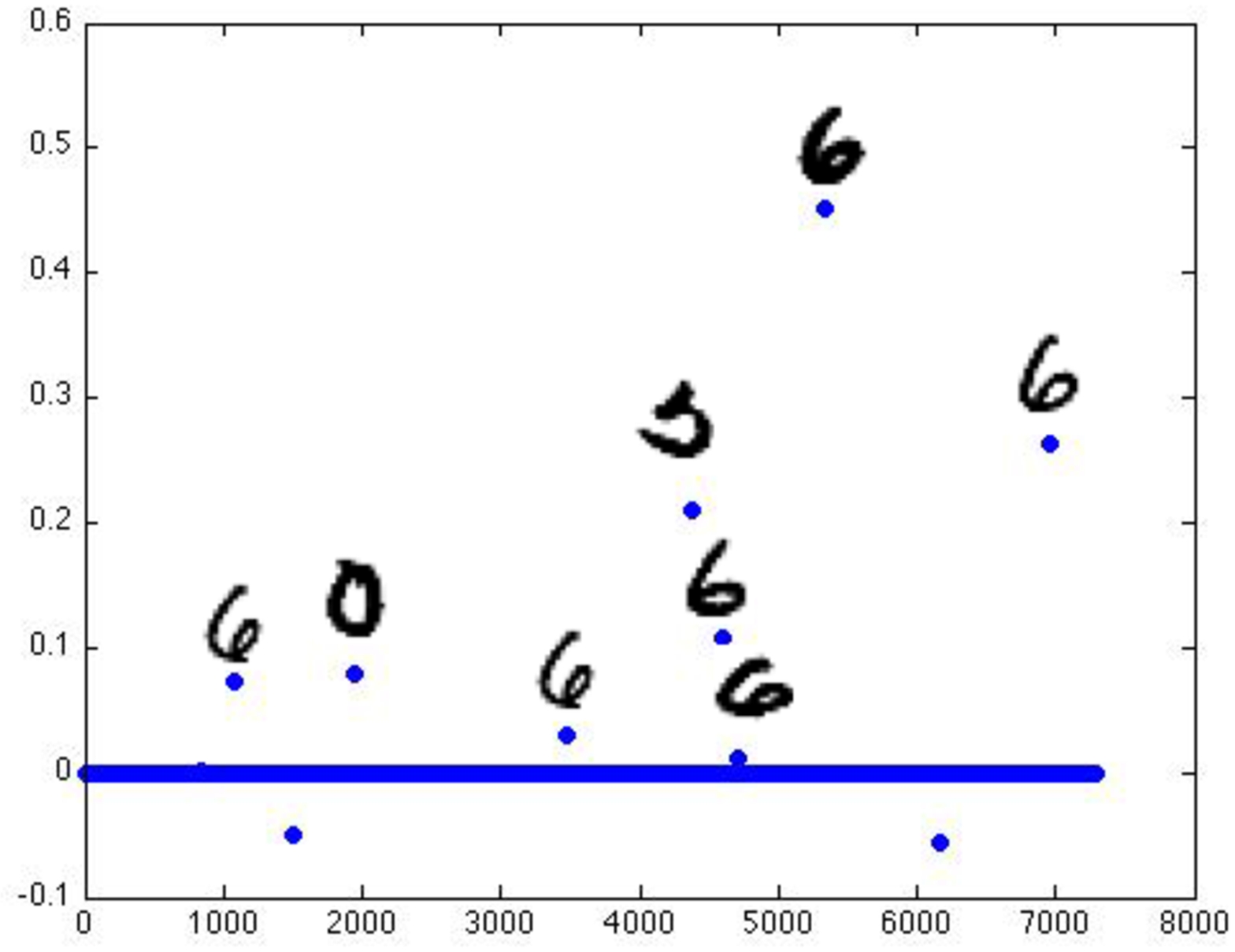}%
\end{center}
\caption{Coefficients of  $\Tilde{\Tilde{\theta}}^{\textsc{es}}_T$, $T=10,000$ and the corresponding image.} \label{FIG:multi}
\end{figure}

The aim of this subsection is to illustrate the performance of the
{\sc es} algorithm on a real dataset and to compare it with the
state-of-the-art procedure in sparse estimation, namely the Lasso.
While sparse estimation is the object of many recent statistical
studies, it is still hard to find a freely available benchmark
dataset where $M \gg n$. We propose the following real dataset
originally introduced in \citet{LeCBosDen90} and, in the particular
instance of this paper, obtained from the webpage of the book
by~\citet{HasTibFri01}. We observe a grayscale image of size
$16\times 16$ pixels of the handwritten digit ``6" (see
Figure~\ref{FIG:6_05}) which is artificially corrupted by a Gaussian
noise. Formally, we can write
\begin{equation}
\label{EQ:mod6}
\bY=\mu+\sigma\xi\,,
\end{equation}
where $\bY \in \R^{256}$ is the observed image, $\mu \in
[0,1]^{256}$ is the true image, $\sigma>0$ and $\xi \in \R^{256}$ is
a standard Gaussian vector. Therefore the number of observations is
equal to the number of pixels: $n=256$. The goal is to reconstruct
$\mu$ using linear combinations of vectors $x_1, \ldots, x_M \in
[0,1]^{256}$ that form a dictionary of size $M=7,290$. Each vector
$x_j$ is a $16\times 16$ grayscale image of a handwritten digit from
$0$ to $9$. As a result, $x_j$'s are strongly correlated as
illustrated by the correlation matrix displayed in
Figure~\ref{FIG:cor_digits}. The digit ``6" is a notably hard
instance due to its similarity with the digits ``0" and with some
instances of the digit ``5" (See Figure~\ref{FIG:multi}). Given an
estimator $\hat \theta$, the performance is measured by the
prediction error $|\mu -\bX\hat \theta|_2^2$, where $\bX$ is the $n
\times M$ design matrix formed by horizontal concatenation of the
column vectors $x_1, \ldots, x_M \in \R^n$.

\begin{figure}[h] \centering
\subfigure[True]{%
\includegraphics[width=0.2\textwidth]{6true.eps}}%
\subfigure[Noisy]{%
\includegraphics[width=0.2\textwidth]{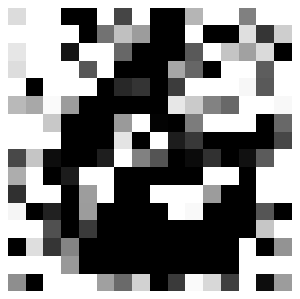}}%
\subfigure[{\sc es}]{%
\includegraphics[width=0.2\textwidth]{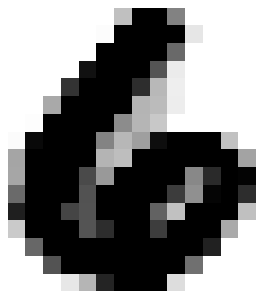}}%
\subfigure[Lasso]{%
\includegraphics[width=0.2\textwidth]{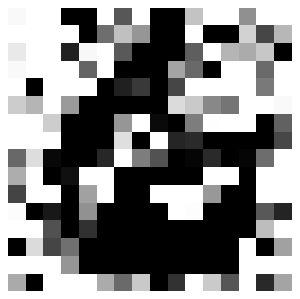}}%
\subfigure[Lasso-Gauss]{%
\includegraphics[width=0.2\textwidth]{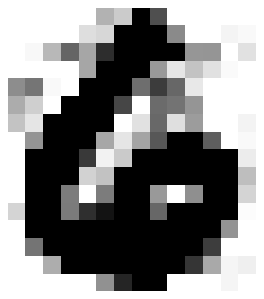}}%
\caption{Reconstruction of the digit ``6" with $\sigma=1.0$.} \label{FIG:6_10}
\end{figure}

Figures~\ref{FIG:6_05} and~\ref{FIG:6_10} illustrate the
reconstruction of this digit by the {\sc es}, Lasso and Lasso-Gauss
estimators for $\sigma=0.5$ and $\sigma=1$ respectively.  The latter
two estimators were computed with fixed regularization parameter
equal to $\sigma\sqrt{8 (\log M)/n}$ and the threshold for the
Lasso-Gauss estimator was taken equal to $\sigma\sqrt{2 (\log M)/n}$.
It is clear from those figures that the Lasso estimator reconstructs
the noisy image and not the true one indicating that the
regularization parameter $\sigma\sqrt{8 (\log M)/n}$ may be too
small for this problem.

\begin{figure}[h]
\psfrag{a}[tr][][.8][45]{{\sc es}}
\psfrag{b}[tr][][.8][45]{Lasso}
\psfrag{c}[tr][][.8][45]{Lasso-G}
\begin{center}
\includegraphics[width=0.5\textwidth]{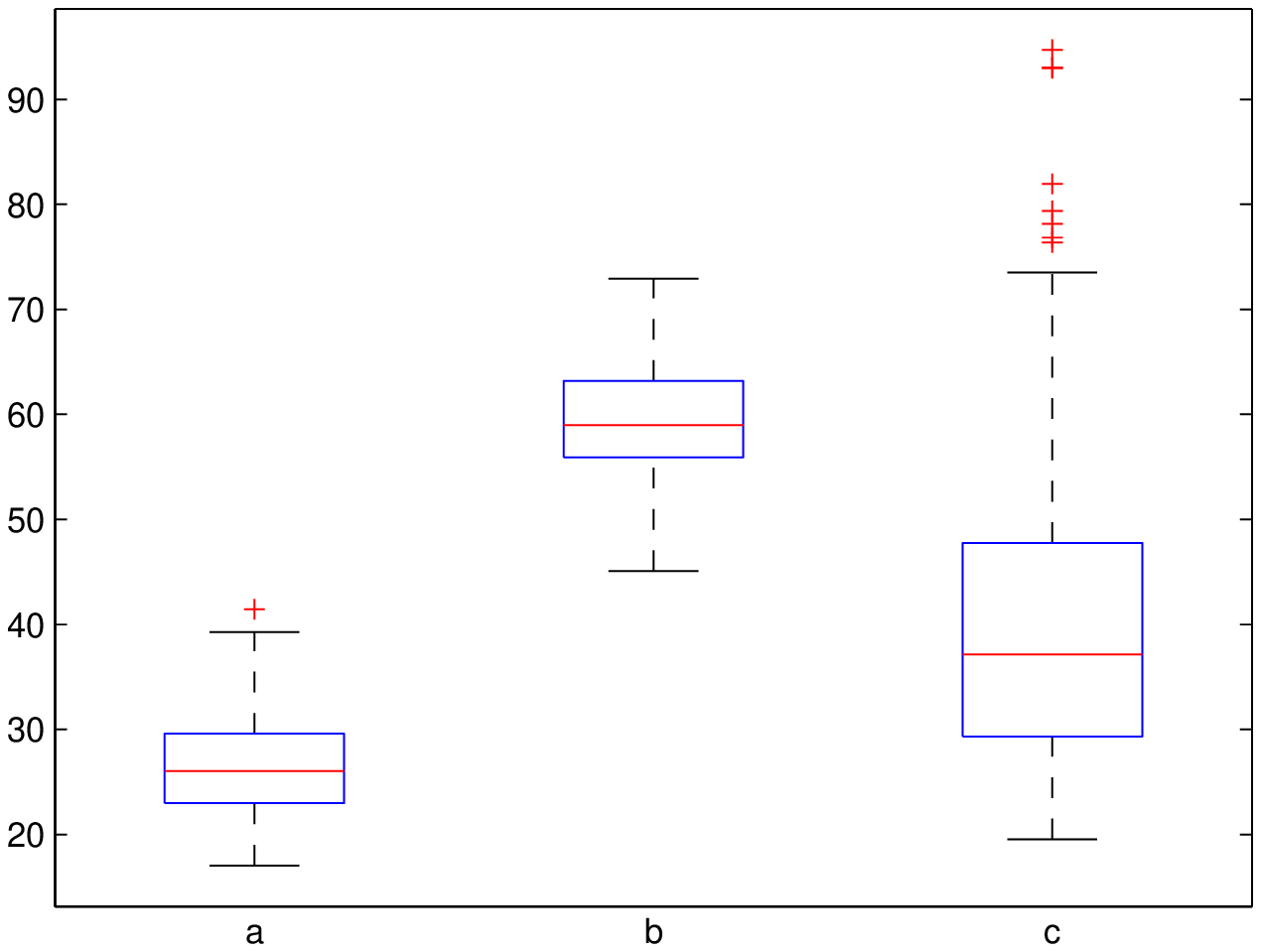}%
\includegraphics[width=0.5\textwidth]{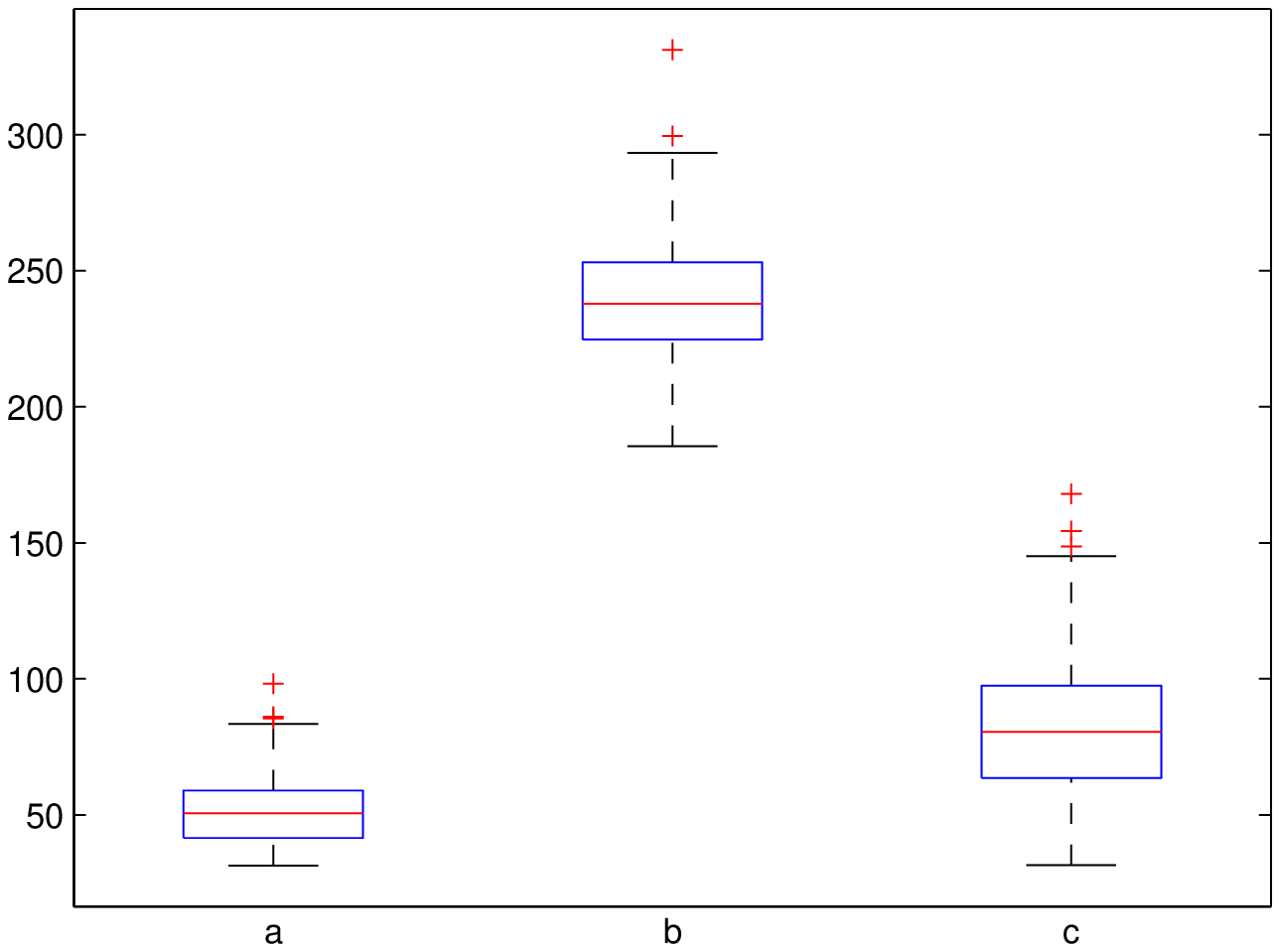}%
\end{center}
\caption{Boxplots of the predictive performance $|\mu -\bX\hat \theta|_2^2$ of the {\sc es}, Lasso and Lasso-Gauss (Lasso-G) estimators computed from 250 replications of the model~\eqref{EQ:mod6} with $\mu$ corresponding to the digit ``6". {\it Left:} $\sigma=0.5$. {\it Right:} $\sigma=1$. Notice that each graph uses a different scale.} \label{FIG:box_digits}
\end{figure}

For both $\sigma=0.5$ and $\sigma=1$, the experiment was repeated
250 times and the predictive performance of {\sc es} was compared
with that of the Lasso and Lasso-Gauss estimators.  The results are
represented in Figure~\ref{FIG:box_digits} and
Table~\ref{TAB:digits}.

\begin{table}[h]
\begin{center}
\begin{tabular}{c|l|l|l}
 & {\sc es} & Lasso & Lasso-Gauss\\
\hline
\hline
$\sigma=0.5$ & ${\bf 26.57}$  &  59.49  &  40.55\\
& {\footnotesize { (4.57)}} &    {\footnotesize (5.28)}  &  {\footnotesize (14.58)  } \\
\hline
$\sigma=1.0$ & ${\bf 51.70}$  &  239.39  &  82.95\\
& {\footnotesize (12.32)} &    {\footnotesize (22.12)}  &  {\footnotesize (24.40)  }\\
\end{tabular}
\end{center}
\caption{Means and standard deviations for $|\mu -\bX\hat
\theta|_2^2$ over 250 realizations of the {\sc es}, Lasso and
Lasso-Gauss estimators to reconstruct the digit
``6".}\label{TAB:digits}
\end{table}

To conclude, we mention a byproduct of this simulation study. The
coefficients of ${\tilde \theta^{\textsc{es}}}$ can be used to
perform multi-class classification following the idea
of~\citet{WriYanGan09}. The procedure consists in performing a
majority vote on the features $x_j$ that are positively weighted by
${\tilde \theta^{\textsc{es}}}$, i.e., such that ${\tilde
\theta^{\textsc{es}}}_j> 0$. For the particular instance illustrated
in Figure~\ref{FIG:6_05}~(c), we see in Figure~\ref{FIG:multi} that
only a few features $x_j$ receive a large positive weight and that a
majority of those correspond to the digit "6".

\section{Appendix}
\label{SEC:appendix}
\setcounter{equation}{0}

\subsection{Lemmas for the upper bound}

The following lemma is obtained by a variant of the ``Maurey
argument", cf. also \citet{Bar93, BunTsyWeg07c, BicRitTsy08} for
similar but somewhat different results.
\begin{lem}
\label{LEM:maurey0} For any $\theta^*\in \R^M\setminus \{0\}$, any
integer $k\ge 1$, and any function $f$ we have
\begin{equation*}
\label{EQ:lem2} \min_{\theta:\begin{subarray}{l}
|\theta|_1=|\theta^*|_1 \\
M(\theta)\le k
\end{subarray}}
\|f-\ff_{ \theta}\|^2\le
\|f-\ff_{\theta^*}\|^2+\frac{|\theta^*|_1^2}{\min(k,
M(\theta^*))}\,.
\end{equation*}
\end{lem}
{\sc Proof.} Fix $\theta^* \in \R^M\setminus \{0\}$  and an integer
$k\ge 1$. Set $K=\min(k, M(\theta^*))$. Consider the multinomial
parameter $p=(p_1, \ldots, p_{M})^\top \in [0,1]^M$ with
$p_j=|\theta^*_j|/|\theta^*|_1, j=1, \ldots, M$, where $\theta^*_j$
are the components of $\theta^*$. Let $\kappa=(\kappa_1, \ldots,
\kappa_M)^\top\in \{0, 1, \ldots, M\}^M$ be the random vector with
multinomial distribution $\cM(K, p)$, i.e., let
$\kappa_j=\sum_{s=1}^{K} \1(I_s=j)$ where $I_1, \ldots, I_{K}$ are
i.i.d. random variables taking value $j \in \{1, \ldots, M\}$ with
probability $p_j$, $j=1, \ldots, M$. In particular, we have
$\cE(\kappa_j)=Kp_j, j=1, \ldots, M$, where $\cE$ denotes the
expectation with respect to the multinomial distribution. As a
result, for the random vector $\bar \theta \in \R^M$ with the
components $\bar \theta_j=\kappa_j\sign(\theta^*_j)|\theta^*|_1/K$
we have $\cE(\bar \theta_j)=\theta^*_j$ for $j=1, \ldots, M$ with
the convention that $\sign(0)=0$. Moreover, using the fact that
${\rm Var}(\kappa_j)=Kp_j(1-p_j)$ and ${\rm Cov}(\kappa_j,
\kappa_l)=-np_jp_l$ for $j\neq l$ \citep[see, e.g.,][eq. (A.13.15),
p.~462]{BicDok06} we find that the covariance matrix of $\bar
\theta$ is given by
$$
\Sigma^*=\cE\left[(\bar \theta - \theta^*)(\bar \theta -
\theta^*)^\top\right]=\frac{|\theta^*|_1}{K}{\rm
diag}(|\theta^*_j|)-\frac{1}{K}|\theta^*||\theta^*|^\top\,,
$$
where $|\theta^*|=(|\theta^*_j|, \ldots, |\theta^*_j|)^\top$. Using
a bias-variance decomposition together with the assumption $\max_j\|f_j\|\le
1$, it yields that,  for any function $f$,
$$
\cE\|f-\ff_{\bar
\theta}\|^2=\|f-\ff_{\theta^*}\|^2+\frac{1}{n}\sum_{i=1}^n
F(x_i)^\top \Sigma^* F(x_i)\le
\|f-\ff_{\theta^*}\|^2+\frac{|\theta^*|_1^2}{K}\,,
$$
where $F(x_i)=(f_1(x_i), \ldots, f_M(x_i))^\top$, $i=1, \ldots, n$.
Moreover, since $\bar \theta$ is such that $|\bar
\theta|_1=|\theta^*|_1$ and $M(\bar \theta)\le K$, the lemma
follows.\epr


\begin{lem}
\label{LEM:maurey} Fix $M, n\ge 1$ and assume that $\max_j\|f_j\|\le
1$. For any function $\eta$ and any constant $\nu>0$ we have
\begin{equation}
\label{EQ:lemMaurey} \min_{\theta \in \R^M}\left\{ \|\ff_{
\theta}-\eta\|^2 + \nu^2\frac{
M(\theta)}{n}\log\left(1+\frac{eM }{ M(\theta)\vee
1}\right)\right\} \le \min_{\theta \in \R^M}\left\{ \|\ff_{
\theta}-\eta\|^2 + {\tilde c} \bar \varphi_{n,M}(\theta)\right\}
\end{equation}
where ${\tilde c}=\left(3+\frac{1}{e}\right)$, $\bar
\varphi_{n,M}(0)=0$ and for $\theta \neq 0$, {
\begin{equation} \label{EQ:varphi} \bar \varphi_{n,M}(\theta)=
\left\{\begin{array}{ll} \min\left[\frac{\nu|\theta|_1}{\sqrt{n}}
\sqrt{\log
\left(1+\frac{eM\nu}{|\theta|_1 \sqrt{n}}\right)}, |\theta|_1^2\right], & {\rm if}\ \langle \ff_\theta, \eta \rangle \le \|\ff_\theta\|^2,\\
\frac{\nu|\theta|_1}{\sqrt{n}} \sqrt{\log
\left(1+\frac{eM\nu}{|\theta|_1 \sqrt{n}}\right)} +
\frac{\nu^2\log(1+eM)}{{\tilde c} n}, & {\rm
otherwise}\,.
\end{array}
\right.
\end{equation}}
\end{lem}
{\sc Proof.} Set
$$
A= \min_{\theta \in \R^M}\left\{ \|\ff_{ \theta}-\eta\|^2 +
\nu^2\frac{ M(\theta)}{n}\log\left(1+\frac{eM }{
M(\theta)\vee 1}\right)\right\}.
$$
It suffices to consider $\R^M\setminus \{0\}$ instead of $\R^M$
since $A\le \|\ff_{0}-\eta\|^2 + {\tilde c}\bar
\varphi_{n,M}(0)=\|\eta\|^2$. Fix $\theta^* \in \R^M\setminus \{0\}$
and define
\begin{equation*}
x^*=|\theta^*|_1 \sqrt{n} /\ell, \quad \text{where} \
 \ell =\nu\sqrt{\log\left(1+\frac{eM\nu }{|\theta^*|_1
\sqrt{n}}\right)}\,.
\end{equation*}
Assume first $x^*\le 1$. In this case we have
\begin{equation}
\label{EQ:xle1}
|\theta^*|_1^2   \le \nu\frac{|\theta^*|_1}{\sqrt{n}}\sqrt{\log\left(1+\frac{eM\nu }{|\theta^*|_1
\sqrt{n}}\right)}\,.
\end{equation}
The previous display yields that $\bar \varphi_{n,M}(\theta^*) = |\theta^*|_1^2$. Moreover, if $\langle \ff_{\theta^*}, \eta \rangle \le \|\ff_{\theta^*}\|^2$, it holds
$$
\|\eta\|^2 \le \|\ff_{\theta^*} -\eta \|^2 + \|\ff_{\theta^*}\|^2 \le  \|\ff_{\theta^*} -\eta \|^2 + |\theta^*|_1^2\,.
$$
As a result,
\begin{equation} \label{EQ:eta2}A\le
\|\eta\|^2 \le \bar \varphi_{n,M}(\theta^*) \quad \mbox{if} \
\langle \ff_{\theta^*}, \eta \rangle \le \|\ff_{\theta^*}\|^2 \
\mbox{and} \ x^* \le 1.
\end{equation}
Set $k^*=\lceil x^* \rceil$, i.e., $k^*$ is the minimal integer
greater than or equal to $x^*$. Using the monotonicity of the
mapping $t\mapsto\frac{t}{n}\log\left(1+\frac{eM }{t}\right)$ for
$t>0$, and Lemma~\ref{LEM:maurey0} we get, for any $\theta^* \in
\R^M\setminus \{0\}$ such that $k^* \le M(\theta^*)\,,$
\begin{eqnarray*}
A&\le & \min_{\theta \in \R^M}\left\{ \|\ff_{ \theta}-\eta\|^2 +
\nu^2\frac{M(\theta)}{n} \log\left(1+\frac{eM }{M(\theta)\vee
1}\right)\right\}
\\
&\le & \min_{1\le k \le M(\theta^*)} \min_{\theta:M(\theta)\le k }\left\{ \|\ff_{ \theta}-\eta\|^2 + \nu^2\frac{k}{n}\log\left(1+\frac{eM }{k}\right)\right\}\\
&\le &  \|\ff_{ \theta^*}-\eta\|^2 +\min_{1\le k \le M(\theta^*)}
\left\{  \nu^2\frac{k}{n}\log\left(1+\frac{eM }{k}\right)+
\frac{|\theta^*|_1^2}{k}\right\}\\
&\le &  \|\ff_{ \theta^*}-\eta\|^2 +
\nu^2\frac{k^*}{n}\log\left(1+\frac{eM }{k^*}\right)+
\frac{|\theta^*|_1^2}{k^*}\,.
\end{eqnarray*}
On the other hand, if $\theta^* \in \R^M\setminus \{0\}$ and $k^*>
M(\theta^*)$, we use the simple bound
\begin{eqnarray*}
A&\le &  \|\ff_{ \theta^*}-\eta\|^2 + \nu^2\frac{M(\theta^*)}{n}
\log\left(1+\frac{eM }{M(\theta^*)\vee 1}\right)
\\
&\le &  \|\ff_{ \theta^*}-\eta\|^2 + \nu^2
\frac{k^*}{n}\log\left(1+\frac{eM }{k^*}\right)\,.
\end{eqnarray*}
In view of the last two displays, to conclude the proof it suffices
to show that
\begin{eqnarray}\label{EQ:proofmaurey1}
\nu^2\frac{k^*}{n}\log\left(1+\frac{eM }{k^*}\right) +
\frac{|\theta^*|_1^2}{k^*} \le {\tilde c} \bar
\varphi_{n,M}(\theta^*)
\end{eqnarray}
for all $\theta^*\in\R^M\setminus \{0\}$. Note first that if $x^* \le 1$, then $k^*=1$ and
$$
\nu^2\frac{k^*}{n}\log\left(1+\frac{eM }{k^*}\right) +
\frac{|\theta^*|_1^2}{k^*}\le\frac{\nu^2\log(1+eM)}{n}+
\nu\frac{|\theta^*|_1}{\sqrt{n}}\sqrt{\log\left(1+\frac{eM\nu
}{|\theta^*|_1 \sqrt{n}}\right)}\le  {\tilde c}
\bar\varphi_{n,M}(\theta^*)\,,
$$
where we used~\eqref{EQ:xle1} in the first inequality. Together
with~\eqref{EQ:eta2}, this proves that $A\le
\|\ff_{\theta^*}-\eta\|^2 + {\tilde c}\bar \varphi_{n,M}(\theta^*)$
for all $\theta^*\in\R^M\setminus \{0\}$ such that $x^* \le 1$.
Thus, to complete the proof of the lemma we only need to consider
the case $x^*>1$. For $x^*>1$ we have
$$
\bar \varphi_{n,M}(\theta^*) \ge \nu\frac{|\theta^*|_1}{\sqrt{n}}\sqrt{\log\left(1+\frac{eM\nu }{|\theta^*|_1
\sqrt{n}}\right)}\,.
$$
As a result, we have
\begin{equation}
\label{EQ:proofmaurey0} \frac{|\theta^*|_1^2}{k^*}\le
\frac{|\theta^*|_1\ell}{ \sqrt{n}} \le  \bar \varphi_{n,M}(\theta^*) \,.
\end{equation}
Moreover, it holds $k^*\le 2x^*=2|\theta^*|_1 \sqrt{n} /\ell$ and
since the function $t\mapsto\frac{t}{n}\log\left(1+\frac{eM
}{t}\right)$ is increasing, we obtain
\begin{equation*}
\frac{k^*}{n}\log\left(1+\frac{eM }{k^*}\right) \le
\frac{2|\theta^*|_1}{ \ell \sqrt{n}} \log\left(1+ \frac{eM \ell
}{2|\theta^*|_1 \sqrt{n}}  \right)\,.
\end{equation*}
Thus, for  $\ell\le \nu$ we have
\begin{equation*}
\frac{k^*}{n}\log\left(1+\frac{eM }{k^*}\right) \le
\frac{2|\theta^*|_1}{ \ell \sqrt{n}} \log\left(1+ \frac{eM
\nu}{2|\theta^*|_1 \sqrt{n}}  \right)\le \frac{2|\theta^*|_1\ell}{
\nu^2\sqrt{n}} \le \frac{2}{\nu^2}  \bar \varphi_{n,M}(\theta^*)\,.
\end{equation*}
For $ \ell> \nu$ we use the inequality $\log(1+ab)\le\log(1+a) +
\log b, \forall \ a\ge 0, b\ge1$, to obtain
\begin{eqnarray*}
\frac{k^*}{n}\log\left(1+\frac{eM }{k^*}\right) &\le&
\frac{2|\theta^*|_1}{ \ell \sqrt{n}} \left[\log\left(1+ \frac{eM
\nu}{2|\theta^*|_1 \sqrt{n}}  \right) + \log\left(\frac{\ell}{\nu}\right)\right] \\
&\le& \frac{2|\theta^*|_1}{ \sqrt{n}} \left(\frac{ \ell}{\nu^2 } + \frac{\log(\ell/\nu)}{\ell} \right) \le
\left(2+\frac1{e}\right)\frac{|\theta^*|_1\ell}{ \nu^2\sqrt{n}}\\
&\le& \left(2+\frac1{e}\right)\frac{1}{\nu^2}  \bar \varphi_{n,M}(\theta^*)\,.
\end{eqnarray*}
Thus, in both cases $\frac{k^*}{n}\log\left(1+\frac{eM
}{k^*}\right)\le (2+1/e)\nu^{-2}  \bar \varphi_{n,M}(\theta^*)$.
Combining this with~\eqref{EQ:proofmaurey0} we
get~\eqref{EQ:proofmaurey1}. \epr

\subsection{Proof of Theorem~\ref{TH:DonJohn1}}

Applying the randomization scheme described in Nemirovski~(2000),
p.211, we create from the sample $y_1,\dots,y_n$ satisfying
\eqref{EQ:gauss_seq_mod} two independent subsamples with
``equivalent" sizes $\lceil n(1-1/\log\log n)\rceil$ and $n - \lceil
n(1-1/\log\log n)\rceil$. We use the first subsample to construct
the {\sc es} estimator and the soft thresholding estimator $\hat
\theta^{\textsc{soft}}$, the latter attaining asymptotically the
rate $\psi_{n}^{01}(\theta)$ for all $\theta\in\R^n$. We then use
the second subsample to aggregate them, for example, as described in
Nemirovski~(2000). Then the aggregated estimator denoted by
$\tilde\theta^*$ satisfies, for all $\theta\in\R^n$,
\begin{eqnarray*}
E_{\theta}|\tilde\theta^*-\theta|_2^2 &\le&
\min\left\{E_{\theta}|\hat \theta^{\textsc{soft}}-\theta|_2^2, \,
E_{\theta}|\tilde \theta^{\textsc{es}}-\theta|_2^2 \right\} +
\frac{C\log\log n}{n}\\
&\le& \min (\psi_{n}^{01}(\theta),\psi_{n}^*(\theta))(1+o(1)) +
\frac{C\log\log n}{n}
\end{eqnarray*}
where $C>0$ is an absolute constant and $o(1)\to 0$ as $n\to\infty$
uniformly in $\theta\in\R^n$. Set
$$
\tilde\psi_{n}(\theta) = \min
(\psi_{n}^{01}(\theta),\psi_{n}^*(\theta)) + \frac{C\log\log n}{n}.
$$
Then \eqref{EQ:THDonJohn11} follows immediately. Next,
$\psi_{n}^{01}(\theta)\ge 2(\log n)/n$, so that for all
$\theta\in\R^n$,
$$
\frac{\psi_{n}^{01}(\theta)}{\tilde\psi_{n}(\theta)}\ge
\frac{\psi_{n}^{01}(\theta)}{\psi_{n}^{01}(\theta) + C(\log\log
n)/n} \ge \frac{2(\log n)/n}{2(\log n)/n + C(\log\log n)/n}\,,
$$
which implies \eqref{EQ:THDonJohn12}. Finally, to prove
\eqref{EQ:THDonJohn13} it is enough to notice that since
$\psi_{n}^{01}(\theta)\ge 2(\log n)/n$,
$$
\frac{\tilde\psi_{n}(\theta)}{\psi_{n}^{01}(\theta)}\le
\frac{\psi_{n}^*(\theta)+ C(\log\log n)/n}{\psi_{n}^{01}(\theta)}
\le \frac{\psi_{n}^*(\theta)}{\psi_{n}^{01}(\theta)} +
\frac{C\log\log n}{2\log n}
$$
and to use \eqref{EQ:THDonJohn2}.


\subsection{Proof of Theorem~\ref{TH:low}}
\label{sub:proof_low}

Clearly \eqref{low1} follows from \eqref{low2} since in the latter
$\eta$ is fixed and equal to one particular function
$\eta=\ff_\theta$.

We now prove \eqref{low2}. Let $\cH=\{f_1, \ldots, f_M\}$ be any
dictionary in $\cD(S\wedge \bar m, { \kappa})$
with the corresponding $\underline \kappa$ and $\bar \kappa$ such
that $\underline \kappa/\bar \kappa= \kappa$. For any $k\in \{1,
\ldots, M\}$, let $\Omega_k$ be the subset of $\cP=\{0,1\}^M$
defined by
\begin{equation}\label{EQ:omegak}
\Omega_k:=\left\{\pp \in \cP\,:\,  |\pp|=k \right\}.
\end{equation}
We consider the class of functions
\begin{equation*}
\cF_k(\delta):=\left\{f=\frac{\delta}{k}\tau\ff_\pp \ : \ \pp \in
\Omega_{k} \right\}\,,
\end{equation*}
where $0<\tau\le 1$ will be chosen later. Note that functions in
$\cF_k(\delta)$ are of the form $\ff_\theta$ with
$\theta\in\R^M_+\setminus\{0\}$, $M(\theta)=k$ and
$|\theta|_1=\tau\delta \le \delta$. Thus, to prove~\eqref{low2}, it
is sufficient to show that, for any estimator $T_n$,
\begin{equation} \label{EQ:LBFl}
 \sup_{\eta \in \cG}   E_\eta\|{T}_n -\eta\|^2 > c_*\kappa \zeta_{n,M, \rank{\bX}}(S,\delta)\,,
\end{equation}
for some subset $\cG \subset \cF_{\bar S}(\delta)$ where $\bar S=
 [ S\wedge\bar m \wedge (M/2)]$ and $[\cdot]$ denotes
the integer part. Note that $\bar S\ge 1$ since $M\ge
2$ and $ S\wedge\bar m \ge 1$.

In what follows we will use the fact that for $f,g\in \cF_{\bar
S}(\delta)$ the difference $f-g$ is of the form $\ff_\theta$ with
some $\theta\in \cP_{2\bar S}$, so that in view
of~\eqref{EQ:class_dict}, $\|f-g\|^2$ is bracketed by the multiples
of $|\theta|^2_2$ with this value of $\theta$.

We now consider three cases, depending on the value of the
integer $m$ defined in~\eqref{EQ:def_m}.

\underline{Case $(1)$: $m=0$}. Use Lemma~\ref{LEM:BM01} to construct
a subset $\cG_{(1)} \subseteq \cF_1(\delta)\subseteq
\cF_{\bar S}(\delta)$ with cardinality $s_{(1)}\ge (1+eM)^{C_1}$
and such that
\begin{equation}
\label{EQ:dist}
\|f-g\|^2 \ge \frac{\tau^2 \delta^2 \underline{ \kappa} }{2}\,, \quad \forall\ f,g \in \cG_{(1)}, f\neq g\,.
\end{equation}
Since $m=0$, inequality~\eqref{EQ:def_m} is violated for $m=1$, so
that
\begin{equation}
\label{EQ:delta1} \delta^2 < \frac{\sigma^2}{n}\log\left(1+eM\right)
\le \frac{\sigma^2}{nC_1}\log(s_{(1)})\,.
\end{equation}

\underline{Case $(2)$: $m\ge 1,  S\wedge (M/2) \ge
m$}.  Then $\bar m=m=\bar S$ and $m \le
M/2$, so that we have $\min(m, M-m)=m$, and Lemma~\ref{LEM:BM01}
guarantees that there exists $\cG_{(2)} \subseteq \cF_m(\delta)=
\cF_{\bar S}(\delta)$  with cardinality $s_{(2)} \ge (1+eM /m)^{C_1m
}$ and such that
\begin{equation*}
\|f-g\|^2 \ge \frac{\tau^2 \delta^2 \underline{ \kappa} }{4m}\,, \quad \forall\ f,g \in \cG_{(2)}, f\neq g\,.
\end{equation*}
To bound from below the quantity $\delta^2/m$, observe that from the
definition of $m$ we have
\begin{equation}
\label{EQ:delta2-1}
\frac{\delta^2}{m} \ge \frac{\delta
\sigma}{\sqrt n}\sqrt{\log\left(1+\frac{eM }{m}\right)} \ge \frac{\delta\sigma}{\sqrt n}\sqrt{\log\left(1+\frac{eM\sigma }{\delta \sqrt n}\right)}\,.
\end{equation}
The previous two displays yield
\begin{equation}
\label{EQ:dist_m}
\|f-g\|^2 \ge \frac{\tau^2  \underline{ \kappa}\delta \sigma}{4\sqrt n}\sqrt{\log\left(1+\frac{eM\sigma }{\delta \sqrt n}\right)}\,.
\end{equation}
Note that in this case
$$
m+1>\frac{\delta \sqrt n}{\sigma\sqrt{\log\left(1+\frac{eM }{m+1}
\right)}}\,,
$$
so that
\begin{equation}
\label{EQ:delta2} \frac{\delta^2}{m}\le \frac{2\delta^2}{m+1}
<2(m+1) \frac{\sigma^2}{n}\log\left(1+\frac{eM }{m+1}\right) \le
4m\frac{\sigma^2}{n}\log\left(1+\frac{eM }{m}\right) \le
\frac{4\sigma^2}{nC_1}\log(s_{(2)})\,.
\end{equation}

\underline{Case $(3)$: $m\ge 1, S\wedge (M/2) < m$}. Then
$\bar S=[ S\wedge (M/2)] < m$. Moreover, we have
$\min(\bar S, M-\bar S)=\bar S$ and using Lemma~\ref{LEM:BM01}, for
any positive $\bar \delta \le \delta$  we can construct
$\cG_{(3)}\subseteq \cF_{\bar S}(\bar \delta)$ with cardinality
$s_{(3)}\ge \big(1+eM /\bar S\big)^{C_1\bar S }$ and such that
\begin{equation*}
\|f-g\|^2 \ge \frac{\tau^2 \bar\delta^2 \underline{ \kappa}
}{4 \bar S} \,, \quad \forall\ f,g \in \cG_{(3)}, f\neq g\,,
\end{equation*}
Take
$$
\bar \delta:=\sigma\frac{\bar S}{\sqrt n} \sqrt{\log \left(1+\frac{eM}{\bar S}\right)} \le
\sigma\frac{m}{\sqrt n} \sqrt{\log \left(1+\frac{eM}{m}\right)}\le \delta\,,
$$
where, in the last inequality, we used the definition of $m$. Next,
note that $\bar S=[ S\wedge (M/2)] 1  \ge S/4$ since
$M\ge 2$. Then
\begin{equation}
\label{EQ:dist_S} \|f-g\|^2 \ge  \frac{\tau^2 \underline{
\kappa}\sigma^2\bar S}{4n}\log\left(1+\frac{eM }{\bar S}\right)\ge
\frac{\tau^2 \underline{ \kappa}\sigma^2
S}{16n}\log\left(1+\frac{4eM }{S}\right)
\,.
\end{equation}
In addition, we have
\begin{equation}
\label{EQ:delta3} \frac{\bar \delta^2}{\bar S}= \bar
S \frac{\sigma^2}{n}\log \left(1+\frac{eM}{\bar S}\right) \le
\frac{\sigma^2}{nC_1}\log(s_{(3)})\,.
\end{equation}

Since the random variables $\xi_i, i=1, \ldots, n$ are i.i.d. Gaussian $\cN(0, \sigma^2)$, for any $f, g \in \cG_{(j)}$, $j\in \{1,2,3\}$, the Kullback-Leibler divergence $\cK(P_f, P_g)$ between $P_f$ and $P_g$ is given by
$$
\cK(P_f, P_g)=\frac{n}{2\sigma^2}\|f-g\|^2\le \frac{n\tau^2 \delta_{(j)}^2\bar \kappa }{\sigma^2k_{(j)}} \,,
$$
where $\delta_{(1)}=\delta_{(2)}=\delta, \delta_{(3)}=\bar \delta,
k_{(1)}=1, k_{(2)}=m, k_{(3)}=\bar S$. Using
respectively~\eqref{EQ:delta1} in case $(1)$, \eqref{EQ:delta2}  in
case $(2)$ and \eqref{EQ:delta3}  in case $(3)$, and choosing
$\tau^2=\min(C_1/(32\bar \kappa),1)$ (note that we need $\tau\le 1$
by construction) we obtain
\begin{equation}
\label{EQ:KL} \cK(P_f, P_g)\le \frac{4\tau^2 \bar \kappa}{C_1}\log
s_{(j)} \le \frac{\log s_{(j)}}{8}\,, \quad \forall\ f, g \in
\cG_{(j)}, \ j=1,2,3.
\end{equation}
Combining~\eqref{EQ:dist}, \eqref{EQ:dist_m} and~\eqref{EQ:dist_S}
together with~\eqref{EQ:KL}, we find that the conditions of
Theorem~2.7 in~\citet{Tsy09} are satisfied and use it to
obtain~\eqref{EQ:LBFl}. \epr

\subsection{A lemma for minimax lower bound}

Here we give a result related to subset extraction, which is a
generalization of the Varshamov-Gilbert lemma used to prove minimax
lower bounds (see, e.g., a recent survey in \citet{Tsy09}[Chap. 2]).
For any $M \ge 1$, $k\in \{1, \ldots, M-1\}$, let $\Omega_k^M$ be the subset of
$\{0,1\}^M$ defined by:
$$
\Omega_k^M:=\left\{\omega \in \{0,1\}^M\,:\, \sum_{j=1}^M \omega_j =k\right\}
$$

The next lemma is a modification of \citet[Lemma~4]{BirMas01}. The
difference is that we cover any $M\ge2$, $1\le k\le M$. The result
of \citet{BirMas01} is proved for even integer $k$ such that $M \ge
3k\ge 6$. The price we pay for considering general $M,k$ is only in
terms of constants, which is sufficient for our purposes.

\begin{lem}
\label{LEM:BM01} Let $M \ge 2$ and $1\le k\le M$ be two integers and
define $\bar k=\min(k, M-k)$. Then there exists a subset $\Omega$ of
$\Omega_{k}^M$ such that the Hamming distance $\rho(\omega,
\omega')=\sum_{j=1}^{M}\1(\omega_j \neq \omega_j')$ satisfies
$$
\rho(\omega, \omega')\ge \frac{\bar k+1}{4}, \qquad \forall\,
\omega,\omega' \in \Omega: \ \omega\neq \omega'\,,
$$
and $s=\card(\Omega)$ satisfies
\begin{equation*}
\log(s)\ge C_1\bar k\log\left(1+\frac{eM}{\bar k}\right)\,,
\end{equation*}
for some numerical constant $C_1\ge 9\cdot10^{-4}$.
\end{lem}
{\sc Proof.} (i) Consider first the case where $k=2p$ for some
integer $p\ge 1$ and $M \ge 6p$. Lemma~4 in \citet{BirMas01}
ensures the existence of a subset $\Omega^{(1)}$ of $\Omega_{k}^M$ such that $\rho(\omega, \omega')\ge k/2+1 \ge (k+1)/2$
for any $\omega\neq \omega' \in \Omega^{(1)}$ and
\begin{equation}
\label{EQ:card_Omega_1}\log\left(\card(\Omega^{(1)})\right) \ge p\left[\log(M/p)-\log(16)+1\right]=
\frac{k}{2}\log\left(\frac{eM}{8k}\right)\,.
\end{equation}

(ii) Next, if $k =2p+1$ for some integer $p\ge 18$ and $M \ge 6p+3$,
let $\widetilde\Omega \subset \Omega^{M-1}_{k-1}$ be the set
obtained by Lemma~4 in \citet{BirMas01}. We have $\rho(\omega,
\omega')\ge (k+1)/2$ for any $\omega, \omega' \in \widetilde \Omega,
\omega'\neq \omega$ and
\begin{equation}
\label{EQ:card_Omega_2}
\log\left(\card(\widetilde\Omega)\right) \ge \frac{k-1}{2}\log\left(\frac{e(M-1)}{8(k-1)}\right)\ge\frac{k}{3}\log\left(\frac{eM}{8k}\right)\,,
\end{equation}
where we used the fact that $3\le k\le M$. Define now the set
$$
\Omega^{(2)}=\left\{\omega \in \{0,1\}^M\,:\, \omega=(1,\tilde \omega)\,,\ \tilde \omega \in \widetilde \Omega\right\}\,.
$$
We have $\Omega^{(2)} \subset \Omega^{M}_k$,
$\card(\Omega^{(2)})=\card(\widetilde\Omega)$ and $\rho(\omega,
\omega')\ge (k+1)/2$ for any $\omega, \omega' \in  \Omega^{(2)},
\omega'\neq \omega$.

So far, we have fully covered $M,k$ such that $M \ge 3k$, $k \ge 36$.
We consider now respectively the cases  (iii) $2k\le M<3k, k \ge 72$, (iv) $k\le 71, M\ge 2k$, and (v) $M<2k$.

(iii) If $2k\le M<3k$, $k \ge 72$, let $k'$ be the
integer part of $k/2$: $k'=[k/2] \ge 36$, and observe that $ 3k' \le
M'$ where $M'=M-(k-k') \le M$. Therefore, we can apply the preceding results to ensure that there exists a
subset $\bar \Omega$ of $\Omega_{k'}^{M'}$ such that
$$
\log\left(\card(\bar \Omega)\right)\ge \frac{k'}{3}\log\left(\frac{eM'}{8k'}\right)\,
$$
and $\rho(\omega, \omega')\ge (k'+1)/2$ for any $\omega, \omega' \in
\bar \Omega$, $\omega \neq \omega'$. Since  $k'\ge k/3$, we obtain
\begin{equation}
\label{EQ:card_Omega_3}
\log\left(\card(\bar \Omega)\right)\ge \frac{k}{9}\log\left(\frac{eM'}{8k'}\right)\,.
\end{equation}
To embed $\bar \Omega$ in $\Omega_{k}^M$, define
$$
\Omega^{(3)}=\left\{\omega \in \{0,1\}^M\,:\, \omega=(\underbrace{1,\dots, 1}_{k-k' \ \text{times}}, \bar \omega)\,,\ \bar \omega \in \bar \Omega\right\}\,.
$$
We have $\Omega^{(3)} \subset \Omega^{M}_k$, $\card(\Omega^{(3)})=\card(\bar \Omega)$ and $\rho(\omega,
\omega')\ge (k'+1)/2 \ge (k+1)/4$ for any $\omega, \omega' \in  \Omega^{(3)}, \omega'\neq \omega$.

(iv) If $k\le 71, M\ge 2k$, consider the set $\Omega^{(4)}=\{\omega^{(1)}, \ldots, \omega^{([M/k])}\} \subset \Omega_k^M$, such that, for any $j=1, \ldots,  [M/k]$, the $l$-th coordinate of $\omega^{(j)}$ satisfies $\omega^{(j)}_l=1$ if and only if $(j-1)k+1\le l\le jk$. We have $\rho(\omega,
\omega')=2k\ge (k+1)/4$ for any $\omega, \omega' \in \Omega^{(4)}, \omega'\neq \omega$ and
\begin{eqnarray}
\label{EQ:card_Omega_4}
\log\left(\card(\Omega^{(4)})\right) =\log\left(\left[\frac{M}{k}\right]\right)&\ge& \frac{\log 2}{\log(1+2e)}\log\left(1+\frac{eM}{k}\right) \nonumber\\
&\ge& \frac{k}{71}\frac{\log 2}{\log(1+2e)}\log\left(1+\frac{eM}{k}\right) \nonumber\\
&\ge& 0.005k\log\left(1+\frac{eM}{k}\right).
\end{eqnarray}

Note that (i)--(iv) cover all $M\ge 2k$ and $k\ge1$, and in
these cases $\bar k =k$. We now use~\eqref{EQ:card_Omega_1},
\eqref{EQ:card_Omega_2} and~\eqref{EQ:card_Omega_3} jointly with the
following inequality
\begin{equation*}
\frac{1}{9}\log\left(\frac{x}{8}\right)\ge \frac{\log\left(\frac{3e}{8}\right)}{9\log (1+3e)}\log (1+x)\ge  0.0009\log (1+x)\,,
\quad x \ge 3e\,.
\end{equation*}
This yields the result of the lemma for cases (i), (ii) and (iii)
since in these cases $M/k\ge 3$ and $M'/k' \ge 3$. For case (iv) we
use directly~\eqref{EQ:card_Omega_4}. Thus, the lemma is proved for
$M\ge 2k$.

(v) Finally, if $M<2k$, or equivalently, when $M-k < k$, we can
reproduce all the arguments above with $k$ replaced by $\bar k =
M-k$ which satisfies $2\bar k \le M$. In each case, $i=1, \dots, 4$, we
obtain the subsets $\bar \Omega^{(i)} \subset \Omega_{\bar k}^M$
analogous to $\Omega^{(i)}$ in (i)--(iv). They are uniquely mapped
into $\Omega_{k}^M$ by applying the bijection $\omega \mapsto {\bf
1}-\omega$, where ${\bf 1}=(1\, \ldots, 1)\in \{0, 1\}^M$.
\epr

{\bf Acknowledgement.} We would like to thank Victor Chernozhukov
for a helpful discussion of the paper.
















\end{document}